\documentclass[a4paper,12pt]{article}
\usepackage[utf8]{inputenc}
\usepackage{amssymb}
\usepackage{amsthm}
\usepackage{amsmath}
\usepackage{mathrsfs}
\usepackage{fancyhdr}
\usepackage{nameref}
\usepackage{abraces}
\usepackage[scr]{rsfso}
\usepackage{mathtools}
\usepackage{xfrac}
\usepackage{xcolor}

\title{Symmetry and classification  of solutions to an integral equation in the Heisenberg group  $\HH^n$}

\author{Jyotshana V. Prajapat\\
\textit{\small Department of Mathematics}\\
\textit{\small University of Mumbai, Vidyanagari}\\
\textit{\small Mumbai 400 098, India}\\
\textit{\small jyotshana.prajapat@mathematics.mu.ac.in}\\
\\
  Anoop Skaria Varghese\\
\textit{\small Department of Mathematics}\\
\textit{\small SIES College of  Arts, Science and Commerce}\\
\textit{\small Affiliated to University of Mumbai}\\
\textit{\small Sion, Mumbai 400 022, India}\\
\textit{\small  anoopv@sies.edu.in }}

\date{}

\begin{document}

	\sloppy
	\theoremstyle{definition}
	\theoremstyle{definition}
	\renewcommand{\contentsname}{Table of Contents}
	
	\newcommand{\be} {\begin{equation}}
	\newcommand{\ee} {\end{equation}}
	\newcommand{\Be} {\begin{equation*}}
	\newcommand{\Ee} {\end{equation*}}

	\newcommand{\bn}{\begin{align}}
	\newcommand{\en}{\end{align}}
	\newcommand{\bea} {\begin{eqnarray}}
	\newcommand{\eea} {\end{eqnarray}}
	\newcommand{\Bea} {\begin{eqnarray*}}
		\newcommand{\Eea} {\end{eqnarray*}}
	\newcommand{\pa} {\partial}
	\newcommand{\ov} {\over}
	\newcommand{\al} {\alpha}
	\newcommand{\ba} {\beta}
	\newcommand{\ka} {\kappa}
	\newcommand{\de} {\delta}
	\newcommand{\ga} {\gamma}
	\newcommand{\Ga} {\Gamma}
	\newcommand{\Om} {\Omega}
	\newcommand{\om} {\omega}
	\newcommand{\De} {\Delta}
	\newcommand{\la} {\lambda}
	\newcommand{\si} {\sigma}
	\newcommand{\Si} {\Sigma}
	\newcommand{\La} {\Lambda}
	\newcommand{\no} {\nonumber}
	\newcommand{\noi} {\noindent}
	\newcommand{\na} {\nabla}
	\newcommand{\vp} {\varphi}
	\newcommand{\var} {\varepsilon}
	\newcommand{\ml} {{\mathscr L}}
	\newcommand{\mh} {{\mathscr H}}
	\newcommand{\supp}{\text{supp }}
	\newcommand{\CC} {\mathbb{C}}
	\newcommand{\RR} {\mathbb{R}}
	\newcommand{\R} {\mathcal R}
	\newcommand{\NN} {\mathbb{N}}
	\newcommand{\HH} {\mathbb{H}}
	\newcommand{\tb}{\textbf}
	\newcommand{\GG}{\mathcal{G}}
	\newcommand{\XX}{\mathcal{X}}

\newcommand{\bu}{
	\overset{\Delta}{u}}
	\newcommand{\bw}{
		\overset{ \Delta}{w}}
	
	\newcommand{\buo}{
		\overset{\Delta}{u_0}}
	\addtolength{\parskip}{0.2em}
	
	\newtheorem{theorem}{Theorem}[section]
	\newtheorem{Pro}[theorem]{Proposition}
	\newtheorem{corollary}[theorem]{Corollary}
	\newtheorem{lemma}[theorem]{Lemma}
	\theoremstyle{remark}
	\theoremstyle{definition}
	\newtheorem{remark}[theorem]{Remark}
	\newtheorem{example}[theorem]{Example}
	\newtheorem{definition}[theorem]{Definition}
	
	\numberwithin{equation}{section}

	\def\H{{{I\!\!H}^n}}
	\def\N{{I\!\!N}}
	\def\T{{\mathbb T}}
	\def\C{\mathbb{C}}
	\def\I{\mathscr{ I}}

	\newcommand{\surface}{T}
	

	\def\sqr#1#2{{\vbox{\hrule height.#2pt
				\hbox{\vrule width.#2pt height#1pt \kern#1pt
					\vrule width.#2pt}
				\hrule height.#2pt}}}
	\def\square{\sqr74}
	\def\qed{{\unskip\nobreak\hfil\penalty50\hskip1em
			\hbox{}\nobreak\hfil\square \parfillskip=0pt
			\finalhyphendemerits=0 \par\goodbreak \vskip8mm}}

	\def\Xint#1{\mathchoice
		{\XXint\displaystyle\textstyle{#1}}%
		{\XXint\textstyle\scriptstyle{#1}}%
		{\XXint\scriptstyle\scriptscriptstyle{#1}}%
		{\XXint\scriptscriptstyle\scriptscriptstyle{#1}}%
		\!\int}
	\def\XXint#1#2#3{{\setbox0=\hbox{$#1{#2#3}{\int}$}
			\vcenter{\hbox{$#2#3$}}\kern-.5\wd0}}
	\def\ddashint{\Xint=}
	\def\dashint{\Xint-}
	
	\renewcommand{\theequation}{\thesection.\arabic{equation}}
	\parindent=0mm
	
	\maketitle

\begin{abstract}
In this paper we prove symmetry of nonnegative solutions of the integral equation
\[ u (\zeta ) = \int\limits_{\HH^n}  |\zeta^{-1} \xi|^{-(Q-\al)}   u(\xi)^{p} d\xi \quad  1< p \leq  \frac{Q+\al}{Q-\al},\quad  0< \al <Q \]
on the Heisenberg group $\HH^n = \C^n \times \RR$, $Q= 2n +2$ using the moving plane method and the Hardy-Littlewood-Sobolev inequality proved by Frank and Lieb for the Heisenberg group. For $p$ subcritical, i.e.,  $1< p <  \frac{Q+\al}{Q-\al}$ we show nonexistence of positive solution of this integral equation, while for the critical case, $p = \frac{Q+\al}{Q-\al}$ we prove that the solutions are cylindrical and are unique upto Heisenberg translation and suitable scaling of the function    
\[ u_0 (z,t) =  \left( (1+ |z|^2)^2 + t^2 \right)^{- \frac{Q-\al}{4}}  \quad (z,t ) \in \HH^n. \]
As a consequence, we also obtain the symmetry and classification of nonnegative $C^2$ solutions of the equation 
\[ \Delta_\HH u + u^{p} = 0 \quad \mbox{for }  1< p \leq \frac{Q+\al}{Q-\al} \mbox{  in } \HH^n \]
without any partial symmetry assumption on the function $u$. 
\end{abstract} 

	\section{Introduction}

One of the important results proved  by Jerison-Lee (\cite{JL1}, \cite{JL2}) for the  CR geometry is the proof of  the CR Yamabe problem stated as ``given a compact, strictly pseudoconvex CR manifold, find a choice of contact form for which the pseudohermitian scalar curvature is constant".  They had conjectured in \cite{JL1}  that  the only solutions to the Yamabe problem on the CR sphere $(S^{2n+1}, \theta_0)$ with the standard contact form $ \theta_0 :=  \frac i2(\bar \pa - \pa) |z|^2$  are images of $\theta_0$ under the CR automorphisms of the sphere $S^{2n+1}$ induced by the biholomorphisms of the unit ball in $\C^{n+1}$, subsequently proving it in their paper \cite{JL2}. The Yamabe problem corresponds to classifying $u \in L^p(\HH^n)$, $ p = \frac{2Q}{Q-2}$ where $Q=2n+2$, positive  solutions of the differential equation 
\be\label{dif}
 \Delta_\HH u + u^{\frac{Q+2}{Q-2}  } = 0 \mbox{  in } \HH^n \ee  
where  $\HH^n = \C^n \times \RR$ is the Heisenberg group with the left group action 
\be
\xi \eta = (z, t) (p, s) := (z+p,t+s+2 \text{Im}\langle z, p\rangle) \ee
with $\langle z, p \rangle =  \sum\limits_{j=1}^{n}z_j {\bar{p}_j}$  denoting the Hermitian inner product and anisotropic scalar multiplication denoted by 
\be \delta_s \xi \mbox{ or } s \xi =(sz, s^2t) \mbox{ for } s \in \RR. \ee
 Moreover, if we denote the vector fields generating the Lie algebra of $\HH^n$ as 
\be X_k = \frac{\pa}{\pa x_k}+ 2y_k\frac{\pa}{\pa t}, Y_k =  \frac{\pa}{\pa y_k}- 2x_k\frac{\pa}{\pa t}, T=\frac{\pa}{\pa t} \mbox{ for } k=1,2,...,n\ee 
then,  \be  \Delta_\HH  = \sum_{k = 1}^n ( X_k^2 + Y_k^2) \ee is the sub-Laplacian operator on $\HH^n$.

It was natural to look for a proof of classification of solutions of (\ref{dif}) using the PDE approach of the moving plane method, which was successful in the Euclidean geometry beginning with \cite{gnn}, \cite{gnn1981}, \cite{chenli} and many subsequent symmetry results. Motivated by this, in \cite{ij},   Birindelli and the first author had initiated the study of symmetry of solutions of PDE in the Heisenberg group using the moving plane method and obtained nonexistence of positive solutions of 
\be\label{dif2}
 \Delta_\HH u + u^{p} = 0 \quad \mbox{for }  1< p < \frac{Q+2}{Q-2} \mbox{  in } \HH^n.\ee  
The symmetry of solutions of subcritical as well as critical exponent problems in the Heisenberg group using the moving plane method  has been elusive since it relies heavily on the maximum principle and the invariance of the differential operator $\Delta_\HH$ under the isometries of the underlying space.  Unlike the usual reflections,  a ``Heisenberg reflection" 
(see \cite{ij}) 
\be
(x, y, t) \mapsto (y, x, 2\la- t) \ee

with respect to the plane orthogonal  to the $t$ axis, say ${\mathcal H}_\la := \{ (z, t) \in \HH^n : t = \la \}$
leaves the plane invariant but not fixed. Therefore, the sign of the difference of  $u$ evaluated at a point and its reflected point  could not be determined on the boundary of the half space $\{(z, t) \in \HH^n : t \geq \la \}$ to be able to apply the maximum principle.  Hence,  earlier efforts to prove the symmetry results in bounded and unbounded domains always required an assumption of  partial symmetry for the domain and/or the function under consideration. In particular,  \cite{ij}, \cite{ij1}, \cite{garo-v} and 
many subsequent results required that the function has a {\bf cylindrical symmetry\/}, i.e.
\be\label{cyl} u(z, t) = u(|z|, t) \quad \mbox{ for } (z, t) \in \HH^n\ee
 where $|z| = \langle z , \bar z \rangle^\frac 12$, $\bar z = (\bar z_1, \ldots, \bar z_n)$ for  $z = (z_1, \ldots, z_n)$. 

In a series of papers beginning with   \cite{CLO}, \cite{CLO1}, \cite{lu}, \cite{jin-li},  \cite{le} (and many more),  
symmetry results were obtained in $\RR^n$ using  the moving plane method for the solutions of the integral equations  (and systems of integral equations) of the type
\be\label{m0}
u(x) =  \int \frac{1}{|x-y|^{n-\al} }u(y)^\frac{n+\al}{n-\al} dy . \ee
The results in \cite{CLO} were particularly interesting for us due to its relation with the usual Yamabe problem and the fact that the use of maximum principles was replaced by the Hardy-Littlewood-Sobolev inequality (henceforth referred to as the  HLS inequality). 

In this paper, we extend the techniques of  \cite{CLO} relying on the HLS inequality for the Heisenberg group proved by Lieb \cite{lieb}.   Consider the integral equation
\be \label{m1}
 u (\zeta ) = \int\limits_{\HH^n} G_\al(\zeta, \xi)  u(\xi)^{p} d\xi \quad  1< p \leq  \frac{Q+\al}{Q-\al},\quad  0< \al <Q \ee
 where \be  G_\al(\zeta, \xi) = 
 |\zeta^{-1} \xi|^{-(Q-\al)}  \ee 
with $| \cdot |$ denoting  the Heisenberg norm $|(z,t)| = (|z|^4 +t^2)^{1/4}$. 

Along with the HLS inequality, another component required in the proof of symmetry is  ``reflections"  in the  Heisenberg group.  In \cite{ij}, we had defined  the $\HH$-reflection with respect to the plane orthogonal  to the $t$ axis in the Heisenberg group as
\be\label{ij} (x, y, t) \mapsto (y, x,  2\la -t ) \quad \mbox{ for } (x, y, t) \in \RR^n \times \RR^n \times \RR.  \ee We had also listed maps such as 
\be 
(x, y, t) \mapsto (2\la -x, y, -t - 4\la y) \ee  which leaves the sub Laplacian invariant.  
 Consistent efforts  of trying to understand the Heisenberg geometry  has led  us to conclude that in principle there is  only  one ``reflection" which matters, i.e., the reflection with respect to the plane orthogonal  to the $t$ axis,  ${\mathcal H}_\la := \{ (z, t) \in \HH^n : t = \la \}$ defined by 
\be\label{refl} (z, t) \mapsto (\bar z,  2 \la -t ) \mbox{  or }  (z, t) \mapsto (- \bar z,  2 \la -t ).\ee 
The cylindrical symmetry (\ref{cyl}) will follow considering the invariance of the equation with respect to the reflection map in (\ref{refl}) and the rotation map given by
\be\label{rota}
(z, t) \mapsto (e^{i\theta} z, t)  \quad \mbox { for } \theta \in [0,2\pi]^n \ee
where $\theta = (\theta_1, \ldots, \theta_n), e^{i\theta} z:= (e^{i \theta_1}z_1, \ldots, e^{i \theta_n}z_n)$. See Step 3 in Subsection \ref{criti} for more details. 

 Henceforth, we will refer to any map described in (\ref{refl}) or  \be  (z, t) \mapsto (e^{i\theta}\bar z,  2 \la -t ) \ee as { \em $\HH$-reflection} and a function $u$ is said to be {\em $\HH$-symmetric} with respect to the plane $t= \la$ for some  $\la\in \RR$ iff 
\be\label{hsym}
u(z, t) = u(|z|,2\la -t).
\ee
 Note that the map (\ref{ij}) is the $\HH$-reflection given by $(z,t) \mapsto (i \bar z, 2\la-t)$. 

Indeed, one can generate $\HH$-reflections with respect to any ``horizontal plane" obtained by translation of the horizontal plane ${\mathcal H}_0 =\{(z, t) \in \HH^n : t = 0\}$ at the origin. ${\mathcal H}_0 $ is also referred to as  the {\em horizontal space } at origin since it is identified with the horizontal tangent space spanned by the set of vector fields $ \{ X_j(0), \, Y_j(0) : 1 \leq j \leq n \}$ evaluated at the origin. The horizontal space ${\mathcal H}_{\xi_0}$ at the point $\xi_0 = ( z_0, t_0)=  ((x_0)_1,\ldots, (x_0)_n,  (y_0)_1, \ldots, (y_0)_n,  t) \in \RR^{2n+1} $ is obtained by left translation
 $ \tau_{\xi_0} {\mathcal H}_0 = \{ \xi_0\eta: \eta \in {\mathcal H}_0 \} $ of the horizontal space at the origin. Here $\tau_{\xi_0}$ denotes Heisenberg translation by a point $\xi_0$. Note that \bea
&& {\mathcal H}_{\xi_0} =  \mbox{Span}\{ X_j(\xi), Y_j(\xi): 1\leq j \leq n\} \nonumber\\
 &= & \mbox{Span}\{ ( 0, \ldots, 0,  \underset{\small{(j)}}{1}, 0,\ldots,    \underset{\small{(2n+1)}}{2(y_0)_j}),
   (0,\ldots, 0,  \underset{\small{(n+j)}}{1}, 0,\ldots, \underset{\small{(2n+1)}}{-2(x_0)_j} ): 1\leq j \leq n\},\nonumber\\\eea is a hyperplane in $\CC^n \times \RR$ passing through $\xi_0$ with the usual normal vector $ (-2 y_0, 2 x_0, 1)\in \RR^{2n+1}$. 
 For $(z, t) \in \HH^n$ the composition of following operations will define $\HH$-reflection with respect to the plane ${\mathcal H}_{\xi_0}$: 
\be \left. \begin{array}{lll} 
\xi =(z, t) & \mapsto & {\xi_0}^{-1}\xi = (z-z_0, t-t_0 -2 \mbox{Im}[\sum\limits_{j=1}^{n}(z_0)_j {\bar{z}_j}] ) \\
&\mapsto &  (\overline{z-z_0}, -t + t_0 + 2 \mbox{Im}[ \sum\limits_{j=1}^{n}(z_0)_j {\bar{z}_j}]) \\
&& ( \mbox{ reflection of the point }   {\xi_0}^{-1}\xi \mbox{ with respect to the plane } {\mathcal H}_0 )\\
& \mapsto & (\overline{z-z_0} + z_0, -t + 2t_0 + 2 \mbox{Im}[ \sum\limits_{j=1}^{n}(z_0)_j {\bar{z}_j}] + 2 \mbox{Im}[ \sum\limits_{j=1}^{n}(z_0)_j ( z- z_0)_j]).
\end{array} \right\} \label{href}
\ee
Let us denote the $\HH$-reflections with respect to the planes $ {\mathcal H}_{\xi_0}$ as ${\mathcal R}_{\xi_0}$  where 
\be {\mathcal R}_0 (z, t) = ( \bar z, -t) \quad \mbox{ for }   (z, t) \in \HH^n\ee
and by \[{\mathcal R}_{\xi_0,\theta} = \xi_0{\mathcal R}_{0,\theta}(\xi_0^{-1} \xi)\] 
where ${\mathcal R}_{0,\theta} (z, t) = (e^{i\theta} \bar z, -t)$  for $ \theta \in [0,2\pi]^n $ and $  (z, t) \in \HH^n$. \\

We first prove the following symmetry for solutions of (\ref{m1}): 

\begin{theorem}\label{symm}
Let $\si = \frac{Q+\al}{Q-\al}$ , $ 0 < \al < Q$, $1<p\leq \si$ and $u \in L^{p +1}_{loc}(\HH^n)$ be a  nonnegative solution of the integral equation (\ref{m1}). Then\\
(i) for $p = \si$, the limit
\be \lim\limits_{|\xi | \to \infty} | \xi |^{Q-\al} u (\xi )  = u_\infty \mbox{ exists},\quad  0 < \al < Q. \ee 
(ii) there exists ${\xi_0} \in \HH^n$ such that    $u$ is $\HH$-symmetric with respect to the plane ${\mathcal H}_{\xi_0}$ i.e.
\be u( \xi) = u\circ {\mathcal R}_{{\xi_0}, \theta}  (\xi) \mbox{ for all } \theta = (\theta_1, \ldots, \theta_n)\in [0,2\pi]^n \mbox{ and } \xi \in \HH^n.
\ee
In particular, $u$ has {\bf cylindrical symmetry\/} up to a Heisenberg translation.\\    
(ii) for $1< p < \si$,   necessarily,  $u$ is $ \HH$-symmetric with respect to the plane $t = 0$, i.e.,
\be\label{s1} u(z, t) = u( |z|, -t ). \ee 
\end{theorem}
To our knowledge, our result is the first of its kind where the moving plane method has been adapted to a setting of a non commutative group.   We have succeeded in proving complete symmetry of solution of integral equation in $\HH^n$  without assuming any partial symmetry or condition at infinity. The ideas in this paper can  contribute to understanding symmetries in  non commutative geometries.  We plan to extend our result to Carnot group and bounded domains in future works. 

Clearly, the condition of $\HH$-symmetry  with respect to the plane ${\mathcal H}_0$  implies that the solution $u$ of (\ref{m1}) is cylindrical as well as even in the $t$-variable (from (\ref{s1}) ),  i.e.,
\be u(z, t) =  u (|z|, t)  = u( |z|, -t) \mbox{ for } (z,t) \in \HH^n. \ee
However, the invariance of the solution $u$ with respect to the reflections 
\be
{\mathcal R}_{\xi_0, \theta}(z, t) =  (e^{i\theta}(\overline{z-z_0}) + z_0, -t + 2t_0 + 2 \mbox{Im}[ \sum\limits_{j=1}^{n}(z_0)_j {\bar{z}_j}] + 2 \mbox{Im}[ \sum\limits_{j=1}^{n}(z_0)_j e^{- i\theta} ( z- z_0)_j])\ee where the term  $ 2 \mbox{Im}[ \sum\limits_{j=1}^{n}(z_0)_j {\bar{z}_j}] + 2 \mbox{Im}[ \sum\limits_{j=1}^{n}(z_0)_j e^{- i\theta} ( z- z_0)_j]$ can be written as 
\Bea
&= &
2 \sum\limits_{j=1}[(y_0)_i (x-x_0)_i - (x_0)_i(y - y_0)_i ]+ \\
&& 2 \sum\limits_{j=1}\left\{ \cos \theta_i [ (x_0)_i(y - y_0)_i + (y_0)_i (x-x_0)_i] - \sin \theta_i [ (x_0)_i (x-x_0)_i - (y_0)_i(y - y_0)_i ]\right\}\Eea
{\em does not } imply that  $u$ is cylindrical about the point $\xi_0$ as the term in $(2n+1)$-th variable also depends on $z$. But  we still note that the solution $u$ depends in the first $2n$ variables only on the distance of $z$ from $z_0$.
Simply put, $v( \zeta) = u(\xi_0^{-1}\zeta)$ is cylindrical.




Typically, to initiate the moving plane method we require the function to have suitable growth condition at infinity. Hence, we work with the CR type inversion of the function $u$ defined as 
\be\label{bu0} \bu( z,t )  =   \frac{1}{| (z,t) |^{Q-\al}} u \bigg( \frac{z}{ \om}, - \frac{t}{|\om|^2}\bigg) \mbox{ with } \om = t + i |z|^2 \quad \mbox { for } (z,t) \in \HH^n \setminus\{0\}\ee
which also satisfies the equation (\ref{m1}) for $p=\si$ (see Lemma \ref{cri}). Note that $\bu$ is the CR inversion in the case $\al = 2$.   In the case $1 < p < \frac{Q+\al}{Q-\al}$ we see that $\bu$ satisfies the equation
\be 
\bu(\zeta ) = \int\limits_{\HH^n} G_\al(\zeta, \xi) \frac{\bu(\xi)^{p} }{ |\xi|^{(Q+\al) - p(Q-\al)}   }   \,d\xi  \mbox{  in } \HH^n \setminus\{0\} \ee
which is of the form \be  v(\zeta ) = \int\limits_{\HH^n}  G_\al(\zeta, \xi) K (|\xi|) v(\xi)^{p} \, d\xi  ,\ee
 similar to the weighted integral equation studied by Chen-Li-Ou in \cite{CLO1} in $\RR^n$. Extension of their results in \cite{CLO1} for supercritical $p$ and general $K$ as well as a study of singular solutions  will appear in forthcoming paper. Symmetry of solutions of integral equations  in bounded domains in Euclidean space is also well studied (see \cite{bdd-symm} and subsequent papers for systems of integral equations on bounded domains). Similar results for bounded domains in  the Heisenberg group and $\HH$-type groups will appear in \cite{prajapat}.

  Here, for the subcritical case $1<p < \frac{Q+\al}{Q-\al}$ with the special case of  $K(|\xi|) = \frac{1}{|\xi|^{(Q+\al) - p(Q-\al)}   }$, we conclude nonexistence of the positive solution of (\ref{m1}) and hence (\ref{dif2}). 
\begin{theorem} {\bf (Nonexistence) }\label{non-exist}
	If $u\in L^{p +1}_{loc}(\HH^n)$ is a nonnegative solution of (\ref{m1}) with $ 1< p < \frac{Q+\al}{Q-\al}$ then $u \equiv 0$.
\end{theorem}

The equation (\ref{m1}) is invariant under the scaling 
\be\label{sca} u_s(z,t) := s^{\frac{Q-\al}{2}} u(sz, s^2 t ) \quad \mbox{ for } { s > 0.} 
\ee
and group translation. See Lemma \ref{ing} for a proof.  For $p = \frac{Q+\al}{Q-\al}$, it can be verified  that \be\label{std0} u_0 = C_0 | \om + i|^{- \frac{Q-\al}{2}} \ee  is a solution of (\ref{m1}), i.e., $u_0$ solves
\be\label{mcrit} u(\xi) =   \int\limits_{\HH^n} G_\al(\xi, \eta)  u(\eta)^{\frac{Q+\al}{Q-\al}} \, d\eta  \ee
We will henceforth refer to $u_0$ as the standard  solution of (\ref{mcrit}).   We classify positive solution of (\ref{m1}) for $p = \frac{Q+\al}{Q-\al}$ as follows. 
\begin{theorem}\label{uniq-th}{\bf (Uniqueness) } Any positive solution $u$  of (\ref{mcrit}) 
 is obtained by a translation and a scaling of the standard solution $u_0 =  C_0| \om + i|^{- \frac{Q-\al}{2}}$. 
 \end{theorem} 
Our proof of uniqueness is different from \cite{CLO} as our symmetry result gives us that  a solution $u$ of (\ref{mcrit})  is a function of two variables. However, we succeed by appealing to the properties satisfied by the standard solution $u_0$. In the process of proving uniqueness,  we first prove the following inversion symmetry. 
\begin{theorem}\label{tinv} Let $u$ be a cylindrical solution of (\ref{mcrit}).  Then there exists $s > 0$ such that 
\be u(r, t) = \frac{s^{Q-\al}}{\rho^{Q-\al}} u\left(\frac{s^2 r}{\rho^{2}}, \frac{s^4 t}{\rho^{4}}\right) \ee
i.e., $u$ is CR inversion symmetric with respect to the CC sphere $\pa B(0,s)$ of radius $s$.
 \end{theorem} 
Here, by a CC sphere $\pa B(0, s)$ we  mean the Carnot-Caratheodory sphere which is the boundary of the open ball 
\[B(0,s) := \{ (z,t) \in \HH^n : (|z|^4 + t^2)^{1/4} < s\}.\]  
Note that $\bu$ defined in (\ref{bu0}) is the CR type inversion with respect to the unit CC sphere $\pa B(0, 1)$. 
The proof of Theorem \ref{tinv} is a consequence of Proposition \ref{inv-sym} and Corollary \ref{3.1} in Section \ref{unq}. 

Since the fundamental solution for the sub Laplacian $\Delta_\HH$ is $|\xi|^{-(Q-2)}$ (see \cite{folland-funda}), the integral equation (\ref{m1}) with $p =  \frac{Q+2}{Q-2} $ and the differential equation (\ref{dif}) are equivalent under suitable regularity assumptions on $u$. Hence, the  Theorem \ref{uniq-th}  finally gives the classification of solutions of the CR Yamabe problem (\ref{dif}). The Liouville theorem for the subcritical case using the moving plane method for the sub Laplacian $\Delta_\HH$ follows from Theorem \ref{non-exist}, {\em without  the condition of cylindrical symmetry} assumed in \cite{ij}. The results for the differential equations associated with the integral equation (\ref{m1})  can be summarised as follows:

\begin{theorem}\label{Yum} Let $u$ be a nonnegative $C^2$ solution of
\be\label{dif3}  \Delta_\HH u + u^{p} = 0 \quad \mbox{for }  1< p \leq \frac{Q+2}{Q-2} \mbox{  in } \HH^n.\ee  
Then the following holds:\\
(i) Symmetry: Any solution of (\ref{dif3}) satisfies $u(z,t) = u(|z|, -t)$ up to a Heisenberg translation ;\\
(ii) Uniqueness: For $p =  \frac{Q+2}{Q-2} $, any positive solution $u$ of (\ref{dif3}) is the standard solution 
\be  u_0 =  C_0| \om + i|^{- \frac{Q-2}{2}} \ee  up to a translation and suitable scaling;\\ 
(iii)Nonexistence: For  $ 1 <p <   \frac{Q+2}{Q-2}$, the only nonnegative solution of (\ref{dif3}) is $u \equiv 0$.
\end{theorem}

{\bf Remark:\/} In recent paper, \cite{cli} proved the classification of solutions of (\ref{dif3}) with $p =\frac{Q+2}{Q-2}$ i.e.,
\be\label{dif3-}  \Delta_\HH u + u^{ \frac{Q+2}{Q-2}} = 0 \quad  \mbox{  in } \HH^n.\ee  

 for $n =1$ or  for $n \geq 2$ with suitable  condition at infinity. Their proof is based on a classical differential identity of Jerison-Lee (\cite{JL2})  combined with integral estimates.  Precisely, they prove
\begin{theorem}(Catino, Li, Monticelli and Roncoroni) 
Let $u$ be a positive solution to (\ref{dif3-}). Then 
\be\label{ustd}
u(z, t) = U_{\la, \mu}(z, t) = \frac{C}{| t + i |z|^2 + z \cdot \mu + \la|^n} 
\ee for some $\la \in \CC$, $\mu \in \CC^n$ such that $Im(\la ) > \frac{|\mu|^2}{4}$. 
\end{theorem} 
and 
\begin{theorem}(Catino, Li, Monticelli and Roncoroni)  
Let $u$ be a positive solution to (\ref{dif3}) in $\HH^n$, $ n\geq 2$ such that
 \be \label{inf-cond}
 u(\xi) \leq  \frac{C}{1+ |\xi|^{\frac{Q-2}{2}}} \quad \mbox{ for all } \xi \in \HH^n \ee
 for some $C > 0$. Then $u$ is of the form (\ref{ustd}).
\end{theorem}
Recall that here $Q = 2n + 2$. We also mention \cite{jj} where the authors prove Liouville theorem for (\ref{dif3-}) under pointwise conditions or integral conditions at infinity:
\begin{theorem}( Flynn and  V\'etois)
Let $n\geq  2$ and $u$ be a positive solution to (\ref{dif3-}) satisfying 
 \be \label{inf-cond}
u(\xi) \leq  \frac{C}{1+ |\xi|^{\frac{n-2}{2}}} \quad \mbox{ for all } \xi \in \HH^n \setminus \{(0,0)\}. \ee
 Then $u$ is of the form (\ref{ustd}).
\end{theorem} 
 \begin{theorem}( Flynn and  V\'etois)
Let $n\geq  2$ and $u$ be a positive solution to (\ref{dif3-}) such that
\be
 \int\limits_{B_R(0)} u^q \leq  C R^2  \quad \mbox{ for all } R>1 \ee
 for some constants $C > 0$ and $q \in (\frac{2n+1}{ n} , \frac{2n+2 }{n}]$. Then $u$ is of the form (\ref{ustd}).
\end{theorem}
We refer to \cite{cli} and \cite{jj} for more details and interesting use of the integral identities and estimates to classify the solutions of (\ref{dif3-}).  Furthermore, in \cite{maou} the authors had proved (iii) i.e., the non existence of positive solutions of (\ref{dif3}) for $1 <p< \frac{Q+2}{Q-2}$ again using a generalized version of Jerison-Lee identity.

Our proof of the complete classification Theorem \ref{Yum}  does not require any extra condition on the function $u$, or any limitations on the dimension of the space considered. Also, for the subcritical case, our proof is a consequence of the symmetry result using the moving plane method. 

Furthermore, from our proof of Theorem \ref{symm}, we can also conclude the following extensions of the results in \cite{gnn1981} to the Heisenberg group. 
\begin{theorem}
 Let $u$ be a positive $C^2$ solution of
\be\label{dif4} \left.\begin{array}{rcl} 
\Delta_\HH u + u^{ \frac{Q+2}{Q-2}} = 0 &   \mbox{  in } & \HH^n \setminus \{0\} \mbox{ with } \\
u(\xi ) \to \infty &\mbox{ as} & |\xi | \to 0 \mbox{ and }\\
  u(\xi ) = O(|\xi|^{Q-2})   &\mbox{ as} & |\xi | \to \infty.
\end{array} \right\} \ee
Then, $u$ is $\HH$-symmetric with respect to ${\mathcal H}_0$ and is decreasing in the $t$ variable i.e. $u_t < 0$ for $t > 0$. 
\end{theorem} 
{\bf Remark:\/}  The conclusion of the solution $u$ being decreasing in the $t$ variable follows once we have proved the symmetry and then applying the moving plane method to the cylindrical solution as in \cite{ij}. 
\begin{theorem}
 Let $u$ be a positive $C^2$ solution of
\be\label{dif5}  \Delta_\HH u + u^{ \frac{Q+2}{Q-2}} = 0 \quad  \mbox{  in } \HH^n \setminus \{0\}\ee  
with singularities at origin and infinity such that 
\be\left.\begin{array}{rcl}
u(\xi ) \to \infty &\mbox{ as} & |\xi | \to 0\\
|\xi|^{Q-2} u(\xi ) \to \infty &\mbox{ as} & |\xi | \to \infty.\\
\end{array} \right\} \ee
Then, $u$ is $\HH$-symmetric with respect to ${\mathcal H}_0$.
\end{theorem}
The proof of this theorem follows from arguments similar to that of proof of Theorem 4 on page 383 in \cite{gnn1981}. 
Generalizations of above results for positive solutions of the equation
\be 
\Delta_\HH u + g(|\xi|, u(\xi)) = 0    \ee   in  $ \HH^n$ or  $ \HH^n \setminus \{0\}$ with suitable conditions on $g$ 
will be studied in future.  \\

For $0< \al < Q$,  it was shown in \cite{cowling} (see also \cite{dooley}, \cite{thanga}) that  $|\xi|^{-(Q- \al)}$ is the fundamental solution of the conformally invariant fractional powers of the sub Laplacian $\Delta_\HH$, which we will  denote by  ${\mathcal L}_{\frac \al2}$ so that $ {\mathcal L}_1 = \Delta_\HH$ for $\al = 2$.    Following \cite{thanga}, for $0 < \al < Q$ we define  ${\mathcal L}_{\frac \al2}$ as 
\be 
 {\mathcal L}_\frac{\al}{2}f(z,t) =(2\pi)^{-n -1} \int\limits_{-\infty}^\infty  \left( \sum_{k=0}^\infty (2|\la|)^\frac{\al}{2} \frac{\Ga\big( \frac{2k+n}{2} +\frac{2+\al}{4}\big)}{\Ga\big( \frac{2k+n}{2} +\frac{2-\al}{4}\big)} f^\la * _\la \phi_k^\la(z)  \right) e^{-i\la t} |\la|^n d\la
\ee
where $\Gamma$ is the Gamma function, $\phi_k^\la$ are the scaled Laguerre functions of the type $(n-1)$ (see pg 7 of \cite{thanga})  and  $*_\la$ is the $\la$-twisted convolution defined in \cite{thanga} as 
\bea f^\la * _\la \phi_k^\la(z)= \int_\CC^n f^\la(z-z') \phi_k^\la (z') e^{\frac i 2 \text{Im} \langle z, z'\rangle} dz'\eea

The Sobolev space  $W^{\frac \al2,2}(\HH^n)$ denotes the collection of all $L^2$ functions $f$ for which $  {\mathcal L}_{\frac \al2}f \in L^2(\HH^n)$. 
The following lemma gives the integral representation of the operator ${\mathcal L}_s$ :
\begin{lemma}( Lemma 5.1 of \cite{thanga}) 
	Let $n\geq 1$ and $0<s= \frac \al2  <1$. Then, for all  $f\in W^{s,2}(\HH^n)$ \be \langle {\mathcal L}_s f, f \rangle = a_{n,s} \int_{\mathbb{H}^n} \int_{\mathbb{H}^n} \frac{|f(\xi)-f(\eta)|^2}{|\xi^{-1}\eta|^{Q+2s}} d\xi d\eta,\ee where $a_{n,s}$ is a positive constant given by  $a_{n,s}=\frac{2^{n-2+3s}}{\pi^{n+1}}\frac{\Ga\big(\frac{n+1+s}{2}\big)^2}{|\Ga(-s)|}.$
\end{lemma}
In section 3 of \cite{thanga}, it was shown that the fundamental solution of $ {\mathcal L}_{\frac \al2}$ is given by 
{\be g_{\frac{\al}{2}}(\xi)=\frac{2^{n+1-3\frac{\al}{2}} \Ga\big(\frac{2n+2-\al}{4}\big)^2}{\pi^{n+1} \Ga(\frac{\al}{2})}|\xi|^{-(Q-\al)}.\ee }

Hence, a function $u$ with suitable regularity satisfies  the integral equation (\ref{m1}) iff it satisfies the  differential equation
\be\label{diff2} {\mathcal L}_{\frac \al2} u+ u^p = 0,\quad  1 < p \leq  \frac{Q+\al}{Q-\al},\quad  0< \al <Q  \mbox{  in } \HH^n. \ee 
Thus, we conclude the following results from  the Theorems \ref{symm},  \ref{non-exist} and \ref{uniq-th}:
 \begin{theorem} Let $u$ be a nonnegative $C^2$ solution of  (\ref{diff2}). Then\\
(i) Symmetry: Any solution of (\ref{diff2}) satisfies $u(z,t) = u(|z|, -t)$ up to Heisenberg translation;\\
(ii) Uniqueness: For $p =  \frac{Q+\al}{Q-\al} $, any positive solution $u$ of (\ref{dif3}) is the standard solution 
\be  u_0 = C_0 | \om + i|^{- \frac{Q-\al}{2}} \ee  upto a translation and suitable scaling;\\ 
(iii)Nonexistence: For  $ 1 <p <   \frac{Q+\al}{Q-\al}$ the only non negative solution of (\ref{diff2}) is $u \equiv 0$.
\end{theorem}

The relation of the sub Laplacian with the Grushin operator defined by 
\be\label{gru}
{\mathcal G} u := \Delta_x u + (s +1)^2 |x|^{2s} \Delta_y u, ~~~{ s > 0 },~~~ (x, y) \in \RR^m \times \RR^k
\ee where $\Delta_x$ and $\Delta_y$ denotes  the usual Laplacian in $\RR^m$ and $\RR^k$
is well known. For, if a function $u(z, t) = u(|z|, t)$ is cylindrical, then 
\be 
\Delta_\HH u  = {\mathcal G} u  \quad \mbox{ for } s = 1 \mbox{ and } m  \mbox{ an  even positive integer }\ee
as the  term 
\be 4 \pa_t(\sum\limits_{i=1}^n y_i\pa_{x_i} - x_i\pa_{y_i}) u = 0 \mbox{ if } u \mbox{ is radial in the } z \mbox{  variable. } \ee    
This discussion would thus be incomplete if we do not  relate our  results to those for the  the semilinear  equations involving the Grushin operator. Precisely, consider a non negative solution of   
\be \label{gr}
\Delta_x u + (s +1)^2 |x|^{2s} \Delta_y u = u^{\frac{Q+2}{Q-2}} ~~~{ s > 0 },~~~ (x, y) \in \RR^m \times \RR^k
\ee where  $Q= m + k(s +1)$ is the homogeneous dimension of $ \RR^m \times \RR^k = \RR^n $.
In \cite{monti}, the authors analyzed the equation (\ref{gr}) and proved symmetry of positive solutions of  (\ref{gr}) using the moving sphere method.  Furthermore, they proved uniqueness of solutions of (\ref{gr}) for any {\em  $s > 0$\/} in the special case  when $m = k =1$   and for the class of $x$-radial function in the case $m \geq 3$ and $k =1$. 
It can be easily seen that for $s = 1$, $k = 1$ and for $m = 2 l$ even integer, the Grushin operator coincides with the Heisenberg sub Laplacian $\Delta_\HH$ acting on cylindrical solutions, i.e., if $u$ is a  cylindrical function $u (z, t) = u(|z|, t)$.

In \cite{garo-v}  the authors studied the critical exponent problem 
\be\label{garo} {\mathcal L} u = - u^{\frac{Q+2}{Q-2}}, \quad u \geq 0, ~~~u \in \mathring{\mathcal D}^{1,2} \ee
in a stratified nilpotent Lie group $G$, also referred to as a Carnot group, where ${\mathcal L}$ is the sub-Laplacian associated with the stratified structure on $G$ and proved symmetry of the cylindrical solution of (\ref{garo}) for $G$ of Heisenberg type. Here, $Q$ is the homogeneous dimension of $G$ and $\mathring{\mathcal D}^{1,2} $ is the closure of $C^\infty$ functions with compact support with respect to the norm $|| u||_{\mathring{\mathcal D}^{1,2}} = ||u||_{L^2} + || Xu||_{L^2}$ where $Xu$ denotes the horizontal gradient.   Again, it can be seen that  for the cylindrical functions, the sub-Laplacian reduces to the Grushin operator (\ref{gru}) with $ s = 1$ and  $m$, $ k $ any positive integer.  
Thus, our uniqueness result (ii) of Theorem \ref{Yum} implies that 
 \begin{theorem}\label{Grush}
 For $s =1$, $k =1$ and $m $ an even integer,  any solution of (\ref{gr}) is  \be  u_0 = C_0 | \om + i|^{- \frac{Q-2}{2}} \ee  upto a translation and suitable scaling.
 \end{theorem}
The proof of uniqueness given in Section \ref{unq} can be extended for Grushin operators for $s =1$, $k \geq 1$ and $m \geq 3 $ as the fundamental solution of the Grushin operator is well known and the differential equation can be associated with the corresponding integral equation.  Details will appear soon. 

 Also, Liouville type   results have been proved for nonlinear elliptic equations involving the Grushin operator in \cite{yu} where the author   proved the nonexistence of positive  solutions of
\be\label{yu} 
{\mathcal G} u + f(u) = 0  ~~~{ s > 0 },~~~ (x, y) \in \RR^m \times \RR^k
\ee 
where $f$ satisfies the conditions\\
($f1$) $f(t) \in  C^0(\RR, \RR)$ is nondecreasing in $(0, \infty)$;\\
($f2$) $g(t) =\dfrac{ f(t)}{t^{\frac{Q+2}{Q-2}}}$ is nonincreasing in $(0, \infty)$  and $g$ is not a constant, where
$Q = m + (s + 1)k$ is the homogeneous dimension.
The $f$ here includes the special case $u^p$ for $1 < p < \frac{Q+2}{Q-2}$, which corresponds to   (iii) of our  
Theorem \ref{Yum}.


The plan of the paper is as follows. In the following section we set up the notations and prove properties inherited by the solutions of (\ref{m1}) due to the invariance of the integral equation (\ref{m1}) under isometries, $\HH$-reflections and the CR inversion.  We  prove Theorem \ref{symm} in  various subsections of Section 3. First, we show symmetry of solutions of (\ref{mcrit}) in Subsection \ref{criti}. Here we also illustrate how to deduce the symmetry of solution from the invariance under $\HH$-reflections. The symmetry of subcritical case $ p < \frac{Q+\al}{Q-\al}$ and non existence are proved in Subsection \ref{subcriti} and \ref{nonex} respectively.  The uniqueness of  solutions of (\ref{mcrit})  is proved in Section 4.

\section{ Notations and  Preliminary results}


Continuing with the notations fixed in the introduction,  the distance of a generic point $\zeta =(z, t) \in \HH^n$ from the origin 
is defined as  \be\rho:= d (\zeta, 0) = | \zeta| =  (|z|^4 +t^2)^{1/4} = (r^4 + t^2)^{1/4}, \ee
where $r = |z|$  denotes the distance of the point $\zeta$ from the $t$-axis. $| \cdot |$ is a norm on $\HH^n$ and  
hence  the Heisenberg distance between two points $\xi$ and $\zeta$ in $\HH^n$ is given by  \[ d (\zeta, \xi) = | \zeta^{-1} \xi|. \]
For future references, the following expression for the distance between points $\zeta = (z_0, t_0)$ and $\eta = (z, t)$ in terms of their coordinates will be useful,
\be\label{distance}
| \zeta^{-1} \eta|^4  =  | (-z_0,-t_0)(z, t)|^4 = | z - z_0 |^4  + | t- t_0 - 2 \text{Im}(z_0 \bar z)|^2. 
\ee

A Carnot Carath\'{e}odary ball or CC-ball  $B(0, \la) \subset \HH^n$ centered at origin is the set 
\[ B(0,\la) : = \{ (z, t) \in \HH^n:  (|z|^4 +t^2)^{1/4} < \la\} .\]
Moreover, the Haar measure on $\HH^n$ is the Lebesgue measure and measure of a set $A \subseteq \HH^n$ will be denoted by $meas(A)$. 

 \be  L^p(\HH^n) := \{ f: \HH^n \to \RR \mbox{ measurable } : \int\limits_{\HH^n} |f(\xi)|^p \,d\xi < \infty\} \nonumber \ee is equipped with the norm 
$|| f||_p :=  \left(\int\limits_{\HH^n} |f(\xi)|^p \,d\xi\right)^{1/p}$. 
Using the invariance of the integral under group translation and scalar multiplication, one can easily verify the following lemma: 

 \begin{lemma}\label{ing} (Invariance under group operations)\\
(i) Scaling:  Let $s\in \RR^+$ and $u$ be a solution of  (\ref{mcrit}). Then $u_s(\xi) =s^{\frac{Q-\al}{2}} u(s \xi )$   also satisfies 
 \be u_s(\xi) =   \int\limits_{\H} G_\al(\xi, \eta)  u_s(\eta)^{\frac{Q+\al}{Q-\al}} \, d\eta.  \ee
(ii) Group translation: 
	Let $u$ be a solution of (\ref{m1}) and $\xi_0 \in \HH^n$. Then $v(\xi) = u( \xi_0  \xi)$ is also a solution of (\ref{m1}).
\end{lemma}

\proof  Since $u$ satisfies the integral equation (\ref{mcrit}), we have 
\[
 u_s(\xi) =  s^{\frac{Q-\al}{2}}	u(s\xi) 
	= s^{\frac{Q-\al}{2}}\int_{\mathbb{H}^n}  G_\al( s\xi, \xi'). u(\xi')^\si d\xi'.
\]
Substitute $ \xi' = s \eta$, then $d \xi'=s^Q d\eta $  and since   \[ G_\al( s\xi, \xi') = G_\al( s\xi, s \eta) = s^{-(Q- \al)} G_\al( \xi, \eta)\] we get 
\Bea u_s(\xi)  & = &  s^{\frac{Q-\al}{2}} \int_{\mathbb{H}^n}  s^{ -(Q-\al)} G_\al( \xi, \eta) u(s\eta)^\si  s^Q d\eta \\
& = &   \int_{\mathbb{H}^n}  G_\al( \xi, \eta) s^{ \frac{(Q+\al)}{2}}  u(s \eta)^\si   d\eta \\
& = & \int_{\mathbb{H}^n}  G_\al( \xi, \eta)  u_s(\eta)^\si   d\eta  \Eea which completes the proof of (i).

(ii) follows easily substituting $ \xi' = \xi_0 \eta$ and observing that $G_\al( \xi_0 \xi,  \xi_0  \eta) = G_\al(\xi, \eta)$.  

\qed
	
	\begin{lemma}\label{inr}
		For $\la \in \RR$, the  equation (\ref{m1}) is invariant under the transformations
		\bea
		{\mathcal X}_\la: (z, t) & \mapsto & (- \bar z, 2 \la -t ), \label{xref}  \\
	{\mathcal R_\theta}:(z, t) & \mapsto & ( e^{i\theta} z,  t ). \label{rota}\eea 
Combining (\ref{xref}) and (\ref{rota}) we see that (\ref{m1}) is invariant under any $\HH$-reflection given by 
 \be (z, t) \mapsto (e^{i\theta}\bar z,  2 \la -t ) .\ee

	\end{lemma}

	\proof Let  $\zeta:=(\tilde{z},\tilde{t}), \xi:=(z,t)$ and $ \xi':=(z',t')$.  To  prove (\ref{xref}), we see that 
since $u$ satisfies the integral equation (\ref{m1}),
\begin{align*}  u({\mathcal X}_\la(\zeta))  =& \int_{\mathbb{H}^n}  G_\al({\mathcal X}_\la(\zeta), \xi') u(\xi')^p d\xi'\\ 
	=& \int_{\mathbb{H}^n}  |(-\bar{\tilde{z}},2\la-\tilde{t})^{-1}.(z',t')|^{-(Q-\al)} u(\xi')^p d\xi'\\
	=& \int_{\mathbb{H}^n} {\Big| |\bar{\tilde{z}}+z'|^4+|-2\la+\tilde{t}+t'+2\text{Im}\langle \bar{\tilde{z}}, z'\rangle|^2 \Big|^{-\frac{Q-\al}{4}}} u(z',t')^p dz'dt'.\end{align*} 

Substituting $(z',t')=(-\bar{z}, 2\la-t)$ gives $dzdt=dz'dt'$, 
$|\bar{\tilde{z}}+z'|=|-\tilde{z}+z|$
and \Be
|\tilde{t}-2\la+ t'+2\text{Im}\langle \bar{\tilde{z}}, z'\rangle|=	|-\tilde{t}+t+2\text{Im}\langle -\tilde{z}, z\rangle|.
\Ee
Hence, 
\begin{align*}
	u({\mathcal X}_\la(\zeta))  =& \int_{\mathbb{H}^n} {\Big||-\tilde{z}+z|^4 +|-\tilde{t}+t+2\text{Im}\langle -\tilde{z}, z\rangle|^2\Big|^{-\frac{Q-\al}{4}}} u(-\bar{z}, 2\la-t)^p dzdt\\
								=& \int_{\mathbb{H}^n}  G_\al(\zeta, \xi). u\circ {\mathcal X}_\la(\xi)^p d\xi.
\end{align*}
\\

 To check (\ref{rota}), define $v(\xi) = v(z, t) = u(e^{i\theta} z, t)$. Then 
 \[ v(\xi)  =  \int_{\mathbb{H}^n}  G_\al((e^{i\theta} z, t), \eta) u(\eta)^p d\eta.\] Substitute $\eta = (e^{i\theta}z', t') = \zeta$. Then $d\eta = dz'dt'$ and $ G_\al((e^{i\theta} z, t), \eta) = G_\al((e^{i\theta} z, t), (e^{i\theta}z', t'))  = G_\al((z, t), (z', t'))$ and the invariance is verified.

\qed 

Next, we collect all the properties of CR type inversion which will be required for our proofs. 	
\subsection{The CR type inversion}\label{CRinv}	
For $\xi \in \HH^n \setminus \{0\}$, let
\bea\label{k1}  \bu( \xi )  & := &  \frac{1}{| \xi |^{Q-\al}} u ( \hat \xi )  \\
\mbox{ where } \hat\xi &:= & \bigg( \frac{z}{ \om}, - \frac{t}{|\om|^2}\bigg) = \bigg( \frac{z\bar \om}{|\om|^2}, - \frac{t}{|\om|^2} \bigg) \mbox{ with } \om = t + i |z|^2. \eea 
denote the CR type inversion of a function $u$.  Note that \be
|\hat \xi| = \frac{1}{|\xi|}, \,\hat \om = - \frac{t}{|\om|^2} + i \frac{|z|^2}{|\om|^2} = - \frac{1}{ \om}
\mbox{ so that }   \hat{\hat\xi} =( - z ,  t) . \ee 
Following properties of the function $\bu$ can be easily verified:\\
(i) The CR-inversion $\bu$ may be singular at the origin. However, if $u$ is continuous then it can be easily seen that 
\be |\bu(\xi)| \leq \frac{C}{|\xi|^{Q-\al}} \mbox{ for all } |\xi| >>0. \ee  
(ii) $u$ is  a cylindrical function  iff  $\bu$ is  a cylindrical function i.e., 
\be \label{crad}
 u(z,t) = u(r, t) \mbox{ iff } \bu(z,t) = \bu(r,t) =\frac{1}{| \xi |^{Q-\al}} u\bigg(\frac{r}{|\om|}, - \frac{t}{|\om|^2} \bigg) \ee
(iii) Unlike the Kelvin transform, the CR type inversion leaves the CC unit sphere  $\pa B(0, 1) =\{ (z, t) \in \HH^n: |z|^4 +t^2 = 1 \}$ invariant and not fixed, i.e., 
\be \bu( z, t) = u (z \bar\om, -t) \quad \mbox{ for } (z, t) \in \pa B(0,1).   \ee
(iv) If $u \in L^p_{loc}(\HH^n)$ then  $ \bu \in L^p_{loc}(\HH^n \setminus \{0\})$.

\begin{lemma}\label{cri}(Invariance of integral equation  (\ref{m1}) under the CR type inversion) 
 If $u$ is a solution of (\ref{m1})
 then  $\bu$ is a solution of  
\be \label{m3}
\bu(\zeta ) = \int\limits_{\HH^n} \frac{1}{ \rho(\zeta, \xi)^{Q - \al} } \frac{\bu(\xi)^{p} }{ |\xi|^{(Q+\al) - p(Q-\al)}   }   \,d\xi  \quad \mbox{for }  \quad 1< p \leq \si= \frac{Q+\al}{Q-\al} \mbox{  in } \HH^n\setminus \{0\}. \ee
In particular, if $p = \si$ then $u$ is a solution of (\ref{m1}) iff $\bu$ is a solution of (\ref{m1}) on $\HH^n \setminus \{0\}$.
 \end{lemma}
\proof     For $\xi \in \HH^n \setminus \{0\} $, 
\[ \bu(\xi) =  \frac{1}{|\xi|^{Q-\al}}  u(\hat \xi) =   \frac{1}{|\xi|^{Q-\al}}   \int\limits_{\HH^n} G_\al(\hat \xi, \eta)  u(\eta)^\si \, d \eta = \frac{1}{|\xi|^{Q-\al}}   \int\limits_{\HH^n} G_\al(\hat \xi, -\eta)  u(-\eta)^\si \, d \eta  \] where $- \eta = (-1)(z,t) = (-z,t)$ is the scaling by $-1$. 
Since  $ \hat{\hat\eta} =(- z , t)$,  substituting $ u(-\eta) = u(\hat{\hat\eta}  ) = |\hat \eta|^{Q-\al} \bu  (\hat \eta)$ and that $d \eta = \dfrac{1}{|\hat \eta|^{2Q}}
d \hat\eta$, we get  
\be  \bu(\xi)  =  \frac{1}{|\xi|^{Q-\al}}   \int\limits_{\HH^n\setminus\{0\}} G_\al(\hat \xi, -\eta) |\hat \eta|^{Q+\al} \bu(\hat \eta)^\si \,  \frac{1}{|\hat \eta|^{2Q}}
d \hat\eta
 =  \int\limits_{\HH^n\setminus\{0\}}\frac{ G_\al(\hat \xi, -\eta)}{|\xi|^{Q-\al}  |\hat \eta|^{Q-\al}} \bu(\hat \eta)^\si \,
d \hat\eta
. \ee
The claim will be proved once we show that
\be \label{cri1}   \frac{ G_\al(\hat \xi, -\eta)}{|\xi|^{Q-\al}  |\hat \eta|^{Q-\al}} = G_\al (\xi , \hat\eta) \ee 
{ i.e., \be |\xi||\hat{\xi}^{-1} (-\eta)|=|\eta|| \xi^{-1}\hat \eta |. \label{Heisen_id}\ee

Denoting $\xi=(z,t), \eta=(p,s), \om=t+i|z|^2$  and $\om'=s+i|p|^2$. Taking $n=1$ for convenience, after computations, we get 
\be |\xi|^4|\hat{\xi}^{-1}( -\eta)|^4=|\om|^2|\om'|^2+1+6|p|^2|z|^2-4\text{Im}(\om'.z\bar{\om}\bar{p})-4\text{Im}(z\bar{p}).\ee 

Similarly, \be |\eta|^4|\hat \eta^{-1} \xi|^4=|\om'|^2|\om|^2+1+6|z|^2|p|^2-4\text{Im}(\om.p\bar{\om'}(-\bar{z}))-4\text{Im}(p(-\bar{z})).\ee
This proves (\ref{Heisen_id}).

\qed

In the following, the constant $C $ denotes a generic positive constant. 
 
\section{Proof of Theorem \ref{symm}}

Here we adapt the moving plane method for the integral equations of \cite{CLO} to the setting of the Heisenberg group.  Apriori, we do not know the behaviour of $u$ at infinity which is essential to begin the moving plane method. However, due to the properties satisfied by  its CR type inversion 
\bea\label{k-1}  \bu( \xi )  & := &  \frac{1}{| \xi |^{Q-\al}} u ( \hat \xi )  \\
\mbox{ where } \hat\xi &:= & \bigg( \frac{z}{ \om}, - \frac{t}{|\om|^2}\bigg) = \bigg( \frac{z\bar \om}{|\om|^2}, - \frac{t}{|\om|^2} \bigg) \mbox{ with } \om = t + i |z|^2, \eea 
 it suffices to prove that $\bu$ is $\HH$-symmetric. We begin with the observation that there exists $R_0 > 0$ and a constant $C_0 > 0$ such that 
\be \label{k-2}
|\bu( \xi ) | \leq   \frac{C_0}{| \xi |^{Q-\al}} \quad \mbox{for all } |\xi| > R_0. \ee
If $u \in L_{loc}^{p+1}(\HH^n)$ for a given $0 < \var < 1$, there exists $R_1 > R_0$ such that 
\be \int\limits_{|\xi| > R} \bu^{p+1} \, d\xi  <  \int\limits_{|\xi| > R} \left(\frac{C_0}{| \xi |^{Q-\al}}\right)^{p+1} \, d\xi  < \var \quad \mbox{ for all } R \geq R_1. \label{k-3}
\ee 
Also note that  $\bu $ is defined in $\HH^n\setminus\{0\}$ with a possible singularity at the origin. We want to prove that after a translation, the function $\bu$ is $\HH$-symmetric with respect to the origin.

  The symmetry of $\bu$ will be obtained by comparing the value of the function $\bu$ at a point $\xi$ and its value at the reflected point $\xi_\la$ with respect to the plane 
\be {\mathcal H}_\la :=\{(z, t) \in \HH^n : t = \la\}. \ee

For $ \la < 0$, denote
\be\label{k1}
\Si_\la = \{ \xi = (z, t)   \in \HH^n : t \geq \la \} \ee
and \be \label{k2} \xi_\la = (\bar z, 2\la -t) \mbox{ for }  \xi = (z, t) \in \Si_\la \ee 
be the reflected point with respect to the plane $\{ t= \la \}  \subset \HH^n$. 
 For $\zeta \in \HH^n$, denote
\be
\Si_{\zeta, \la}= \{ \zeta \xi: \xi \in \Si_\la \} = \tau_\zeta \Si_\la \ee so that $\Si_{0, \la} = \Si_\la$ when $\zeta $ is the origin. 
A point $\eta \in \Si_{\zeta, \la}$ can be written as $\eta = \zeta \xi$ for $\xi \in \Si_\la$ and it can be verified that the reflection of $\eta $ with respect to the plane $\tau_\zeta {\mathcal H}_\la$ is 
\be
 {\mathcal R}_{\zeta, \la}(\eta) =  {\mathcal R}_{\zeta, \la}(\zeta\xi) = \zeta\xi_\la ~~~\mbox{where } \xi_\la \mbox{ is given by } (\ref{k2}). \ee

Define
\be \label{k3}u_\la( \xi) := u( \xi_\la) \quad \mbox{for }  \xi \in \Si_\la.  \ee

Using the Lemmas \ref{ing}, \ref{inr} and \ref{cri}  proved in the previous section, we will first prove the following important lemma.
\begin{lemma}\label{l2.1} If $u$ be a solution of  (\ref{m1}) then for $\zeta \in \Si_\la$, 
	\be\label{r1}
	u (\zeta_\la ) - u (\zeta )  =   \int\limits_{\Si_\la} |\eta|^{-(Q-\al)} \left(u_\la(   \zeta\eta  )^{p} -  u(\zeta \eta)^{p} \right)d\eta 
 +  \int\limits_{\Si_\la}  | \eta_\la |^{-(Q-\al)} \left(  u_\la(\zeta\eta_\la)^{p} -  u(\zeta \eta_\la)^{p} \right)  d\eta. 
 \ee

\end{lemma}

\proof 
For $\la = 0$,  $\zeta_0 = (x_0, -y_0, - t_0)$ is the reflection of $\zeta = ( x_0, y_0, t_0)$  with respect to the plane $t= 0$. The following relations can be verified:\\
(i)  $(\zeta^{-1})_0 = (\zeta_0)^{-1}$.\\
(ii) $(\zeta\xi)_0 = \zeta_0 \xi_0$.\\ 
(iii)  $\zeta_\la^{-1}\xi_\la  = (\zeta^{-1})_0 (\xi)_0 = ( \zeta^{-1}\xi   )_0$. This gives another verification of $G_\al(\zeta_\la, \xi_\la) = G_\al(\zeta, \xi) $.\\
(iv) $ (\zeta_\la)_0  \neq (\zeta_0)_\la$\\
(v)  If $\xi_\la = \zeta \eta $ then $\xi = \zeta_0 \eta_\la$.\\
Using above relations, we have 
\Bea
u (\zeta_\la ) &  = &  \int\limits_{\HH^n} G_\al(\zeta_\la, \xi)  u(\xi)^{p} d\xi\\
&  =&   \int\limits_{\HH^n} |\zeta_\la^{-1} \xi|^{-(Q-\al)} u(\xi)^{p} d\xi   =   \int\limits_{\HH^n}  | \eta |^{-(Q-\al)} u((\zeta\eta)_\la)^{p} d\eta\\
&&\mbox{ using } \zeta_\la^{-1} \xi = ( \zeta^{-1} \xi_\la)_0 \mbox{ and } |\zeta_\la^{-1} \xi| = |( \zeta^{-1} \xi_\la)_0|  = |\zeta^{-1} \xi_\la|\\
& =& \int\limits_{\Si_\la} |\eta|^{-(Q-\al)} u(   (\zeta\eta)_\la  )^{p} d\eta   +    \int\limits_{\Si_\la}  | \eta_\la |^{-(Q-\al)} u((\zeta\eta_\la)_\la)^{p} d\eta\\
& =& \int\limits_{\Si_\la} |\eta|^{-(Q-\al)} u_\la(   \zeta\eta  )^{p} d\eta   +    \int\limits_{\Si_\la}  | \eta_\la |^{-(Q-\al)} u_\la(\zeta\eta_\la)^{p} d\eta.
\Eea

We then write \Bea
u (\zeta )  & = &  \int\limits_{\HH^n} G_\al(\zeta, \xi)  u(\xi)^{p} d\xi  = \int\limits_{\HH^n} |\zeta^{-1} \xi|^{-(Q-\al)} 
u(\xi)^{p} d\xi = \int\limits_{\HH^n} |\eta|^{-(Q-\al)} 
u(\zeta\eta)^{p} d\eta\\ 
& = &   \int\limits_{\Si_{\la}} |  \eta|^{-(Q-\al) } u(\zeta \eta)^{p} d\eta +   \int\limits_{\Si_{\la}}  | \eta_\la|^{-(Q-\al)} u(\zeta \eta_\la)^{p} d\eta.
\Eea
 Therefore, 
 (\ref{r1}) holds.

\qed

Regarding the integrands on the RHS of (\ref{r1}) we have following observations:\\

(i)  In general, the Heisenberg translation does not translate the plane ${\mathcal H}_\la$ parallel to itself, in fact for any $\zeta = (z_0, t_0)$ as long as $ z_0 \neq 0$, the translated plane $\tau_\zeta {\mathcal H}_\la$ will never be parallel to the original plane ${\mathcal H}_\la$.  Hence, for any such $\zeta$,  $\tau_\zeta \Si_\la \cap \Si_\la^c$ will always be {\em nonempty}. \\

(ii) The set $\Si_\la$ can be written as a disjoint union $ \Si_\la = (\tau_\zeta \Si_\la \cap \Si_\la) \cup (\tau_\zeta \Si_\la)^c \cap \Si_\la$, where $ (\tau_\zeta \Si_\la)^c  $ is the reflection $
 {\mathcal R}_{\zeta, \la} \tau_\zeta \Si_\la $ of $\tau_\zeta \Si_\la$ with respect to the plane $\tau_\zeta {\mathcal H}_\la$. Hence every point in $\Si_\la$ can be expressed as $\zeta \eta$ or $\zeta \eta_\la$ where $\eta \in \Si_\la$, which has been used in the representation (\ref{r1}). \\

(iii) To see that this is the correct representation which will help achieve our goal of proving symmetry, 
suppose that we can show that both the integrands  $\left(u_\la(   \zeta\eta  )^{p} -  u(\zeta \eta)^{p} \right) $ and $ \left(  u_\la(\zeta\eta_\la)^{p} -  u(\zeta \eta_\la)^{p} \right)$ does not change sign, say are non negative for all $\zeta \in \Si_\la$. Then, $u(\zeta_\la) \leq u(\zeta)$ in $\Si_\la$.  If,  for some $\zeta \in \Si_\la$,  $u(\zeta) = u(\zeta_\la)$ i.e.,
\[0 =   \int\limits_{\Si_\la} |\eta|^{-(Q-\al)} \left(u_\la(   \zeta\eta  )^{p} -  u(\zeta \eta)^{p} \right)d\eta 
 +  \int\limits_{\Si_\la}  | \eta_\la |^{-(Q-\al)} \left(  u_\la(\zeta\eta_\la)^{p} -  u(\zeta \eta_\la)^{p} \right)  d\eta
    \]
 then we can conclude that 
\be\label{ksym}  u_\la(\zeta\eta)   =    u(\zeta\eta) \mbox{ and }   u_\la(\zeta\eta_\la) = u(\zeta \eta_\la) \quad \mbox{ for all } \eta \in \Si_\la,\ee i.e., $u $ is $\HH$ symmetric in $\tau_\zeta \HH^n$ with respect to the plane $\zeta^{-1}{\mathcal H}_\la$.  

Since we will first prove symmetry of the CR inversion $\bu$, an immediate consequence of the above lemma is 
\begin{corollary} The CR type inversion $\bu$ solves (\ref{m1})  with $p = \si$ for $\zeta \in \Si_\la \setminus \{0\} = \Si_\la^*$ and hence satisfies ,
	\bea\label{i1}
&&	\bu(\zeta_\la ) - \bu (\zeta)  \nonumber\\
&= & \int\limits_{\Si_\la } |\eta|^{-(Q-\al)} \left(\bu_\la(   \zeta\eta  )^{p} - \bu(\zeta \eta)^{p} \right)d\eta 
 +  \int\limits_{\Si_\la }  | \eta_\la |^{-(Q-\al)} \left( \bu_\la(\zeta\eta_\la)^{p} -  \bu(\zeta \eta_\la)^{p} \right)  d\eta. \nonumber
 \\
\eea 
In general, for $ p \leq \si$,  define
	\be 
	f(\xi) := \frac{\bu(\xi)^{p} }{ |\xi|^{(Q+\al) - p(Q-\al)}} 
	\ee
so that (\ref{m3}) becomes 
\be\label{km3}     \bu(\zeta ) = \int\limits_{\HH^n} \frac{1}{ \rho(\zeta, \xi)^{Q - \al} } f(\xi)   \,d\xi  \quad \mbox{for }  \quad 1< p \leq \si= \frac{Q+\al}{Q-\al} \mbox{  in } \HH^n\setminus \{0\}. \ee
Then 
\be\label{i2}
\bu(\zeta_\la ) - \bu (\zeta)  = \int\limits_{\Si_\la } |\eta|^{-(Q-\al)} \left(f_\la(   \zeta\eta  ) - f(\zeta \eta) \right)d\eta 
 +  \int\limits_{\Si_\la }  | \eta_\la |^{-(Q-\al)} \left( f_\la(\zeta\eta_\la) -  f(\zeta \eta_\la) \right)  d\eta
\ee
\end{corollary}

Theorem \ref{symm} will follow from the results proved in the following subsections. We will first prove the symmetry of $\bu$ solution of (\ref{km3}) considering the critical and subcritical case separately. 

\subsection{ Invariance  of solutions of (\ref{km3})  under $\HH$-reflection}\label{criti} 
{\bf 	Case (i)\/} $ {\bf p = \si=  \frac{Q+\al}{Q-\al}}$: 
From (\ref{k-2}),  in a deleted neighbourhood of $0$ we always have that 
\[ \bu(\xi ) \geq \bu(\xi_\la) \mbox{  for all large }-\la>0. \] We thus claim that 
 $\bu(\xi) \geq \bu_\la(\xi) $ in $\Si_\la$ for  all $\la << 0$.  Let \be
E_\la:= \{ \xi \in \Si_\la^*: \bu_\la(\xi) > \bu(\xi)\} \ee
 and for $\zeta \in \Si_\la^*$, denote 
\be\label{seta}
A_{\la, \zeta} := \{\xi \in \Si_\la : \zeta \xi \in E_\la \} \mbox{ and  } 
B_{\la, \zeta} := \{\xi \in \Si_\la : \zeta \xi_\la \in E_\la  \}. \ee
Then 
\be\label{sets}  E_\la = \cup_{\zeta \in \Si_\la}( \tau_\zeta A_{\la, \zeta} \cup \tau_\zeta (B_{\la, \zeta})_\la)\ee  where $(B_{\la, \zeta})_\la$ is reflection of the set $B_{\la, \zeta}$ with respect to the plane ${\mathcal H}_\la$.  
Note that if $\zeta \in E_\la$ then $0 \in A_{\la, \zeta}$. Also, $\zeta^{-1} \notin A_{\la, \zeta} \cup B_{\la, \zeta}$. 

We will prove our result  in the following steps. \\

 \noi{\bf Step 1: There exists $M_0 > 0$ such that for all $\zeta \in \Si_\la^*$,
 	\be   \bu_\la(\zeta\xi) \leq  \bu(\zeta \xi)  \mbox{  and }  \bu_\la(\zeta\xi_\la) \leq \bu (\zeta\xi_\la )
\mbox{ for  }  \xi \in \Si_{ \la} 
  \mbox{ and for  all } \la < -M_0 . \ee
In particular, since $0 \in \Si_\la$,  
	\be   \bu(\zeta_\la) \leq  \bu(\zeta )\mbox{ for  }  \xi \in \Si_{ \la} 
  \mbox{ and for  all } \la < -M_0 . \ee}
\\
The integral (\ref{r1}) can be written as 
\bea
&& \bu (\zeta_\la ) -\bu (\zeta ) \nonumber \\ &\leq &    \int\limits_{  A_{\la, \zeta}         } |\eta|^{-(Q-\al)} \left(\bu_\la (   \zeta\eta  )^{\si} -  \bu(\zeta \eta)^{\si} \right)d\eta 
 +  \int\limits_{B_{\la, \zeta}}  | \eta_\la |^{-(Q-\al)} \left(  \bu_\la (\zeta\eta_\la)^{\si} -  \bu(\zeta \eta_\la)^{\si} \right)  d\eta\nonumber \\
&=&  \int\limits_{ \zeta  A_{\la, \zeta}} |\zeta^{-1}\xi|^{-(Q-\al)} \left(\bu_\la (   \xi  )^{\si} -  \bu( \xi)^{\si} \right)d\xi 
+  \int\limits_{\zeta (B_{\la, \zeta})_\la}  |\zeta^{-1}\xi |^{-(Q-\al)} \left(  \bu_\la (\xi)^{\si} -  \bu(\xi)^{\si} \right)  d\xi. \nonumber \\
\eea
Define \Bea c(\xi) &=& \left( \frac{ \bu_\la( \xi)^\si - \bu(\xi )  )^\si}{ \bu_\la(\xi) - \bu( \xi   )}\right) (\bu_\la(\xi) - \bu( \xi ) ) \mbox{~~if ~~}\bu_\la(\xi) \neq \bu(\xi ) \\
&=0& \mbox{~~otherwise}. \Eea
Now, for $\xi \in E_\la$
\be \bu_\la(\xi)^\si - \bu(\xi )^\si = \si a(\xi)(\bu_\la( \xi) - \bu(\xi )  )\leq \si\bu_\la(\xi)^{\si -1 } ( \bu_\la(\xi) - \bu (\xi ) ) \ee
where $a(\xi)$ is a real number between $\bu_\la(\xi)$ and $\bu(\xi)$ and  we also use that  $a(\xi) < \bu_\la(\xi)$ for $\xi \in E_\la$. 
Therefore,
\bea
\bu (\zeta_\la ) -\bu (\zeta ) 
& \leq&  \int\limits_{\zeta  A_{\la, \zeta}         } |\zeta^{-1}\xi|^{-(Q-\al)} \si\bu_\la(\xi)^{\si -1 } ( \bu_\la(\xi) - \bu (\xi )) \,d\xi \nonumber\\
&& +  \int\limits_{\zeta (B_{\la, \zeta})_\la}   | \zeta^{-1}\xi|^{-(Q-\al)} \si\bu_\la(\xi)^{\si -1 } ( \bu_\la(\xi) - \bu(\xi ) ) \,d\xi \label{t1}\eea
Step 1 will be proved once we show\\
{\bf Claim: $E_\la$ is a set of measure zero.}\\
The following argument is written for $1< p \leq \si$ to avoid repetition in the sub critical case. For a $q > 1 $ ( $q =\frac{2 Q}{Q-\al}$ in case $p = \si$),  multiplying both sides of (\ref{t1})  by $  \left(\bu(\zeta_\la) - \bu(\zeta)\right)^{q-1}$ and integrating on $E_\la$ we get 
\bea \label{-hls}
 &&\int\limits_{E_\la} \left(\bu(\zeta_\la) - \bu(\zeta)\right)^{q} \,d\zeta \nonumber \\
 &\leq & 
  \int\limits_{E_\la}\int\limits_{\zeta A_{\la, \zeta}}     |\zeta^{-1}\xi|^{-(Q-\al)} \si\bu_\la(\xi)^{\si -1 } ( \bu_\la(\xi) - \bu (\xi )  \left(\bu(\zeta_\la) - \bu(\zeta)\right)^{q-1} \, d\xi \,d\zeta \nonumber\\
&&+\int\limits_{E_\la}\int\limits_{\zeta (B_{\la, \zeta})_\la }  |\zeta^{-1}\xi|^{-(Q-\al)} \si\bu_\la(\xi)^{\si -1 } ( \bu_\la(\xi) - \bu (\xi ) \left(\bu(\zeta_\la) - \bu(\zeta)\right)^{q-1} \, d\xi d\zeta \nonumber \\
& \leq &  \int\limits_{E_\la}\int\limits_{E_\la}  |\zeta^{-1}\xi|^{-(Q-\al)} \si\bu_\la(\xi)^{\si -1 } ( \bu_\la(\xi) - \bu (\xi ) \left(\bu(\zeta_\la) - \bu(\zeta)\right)^{q-1} \,d\zeta \, d\xi  \nonumber\\
&&+\int\limits_{E_\la}\int\limits_{E_\la} |\zeta^{-1}\xi|^{-(Q-\al)} \si\bu_\la(\xi)^{\si -1 } ( \bu_\la(\xi) - \bu (\xi ) \left(\bu(\zeta_\la) - \bu(\zeta)\right)^{q-1} \,d\zeta \, d\xi \nonumber \\
& = & 2 \int\limits_{E_\la}\int\limits_{E_\la} |\zeta^{-1}\xi|^{-(Q-\al)} \si\bu_\la(\xi)^{\si -1 } ( \bu_\la(\xi) - \bu (\xi ) \left(\bu(\zeta_\la) - \bu(\zeta)\right)^{q-1} \,d\zeta \, d\xi. \eea

 Recall  the HLS inequality for the Heisenberg group proved by Frank-Lieb \cite{lieb} (in the notations used therein):  
 \begin{theorem}(HLS inequality)
 	Let $0 < \la < Q = 2n +2$ and $p := \frac{2Q}{2Q-\la}$. Then for any $f$, $g \in L^p(\HH^n)$, 
 	\be\label{hls}
 	\left| \iint_{\HH^n \times \HH^n} \frac{ \overline{f(\xi)} g(\eta)}{|\xi^{-1} \eta|^\la} \, d\xi d\eta \right| \leq \left(\frac{\pi^{n+1}}{2^{n-1}n!} \right)^{\la/Q} \frac{n! \Ga((Q-2)/2)}{\Ga^2((2Q-\la)/4)} \| f\|_p \|g\|_p \ee
 	with equality if and only if \be f(\xi)=cH(\de(a^{-1}\xi)),\quad  g(\xi)=c'H(\de(a^{-1}\xi)) \ee for some $c,c'\in \CC, \de>0$ and $a\in \HH^n$(unless $f\equiv0$ or $g\equiv 0$). Here, $H$ is the function given by $H(z,t)= \Big( (1+|z|^2)^2+t^2\Big)^{-(2Q-\la)/4}$.	 \end{theorem}
  Here $ H$ is the standard solution $u_0$ defined in (\ref{std0}). \\


Applying the HLS inequality (\ref{hls}) to (\ref{-hls}), we get 
\[
\int\limits_{E_\la} \left(\bu_\la(\zeta) - \bu(\zeta)\right)^{q} \,d\zeta  \leq  C \si ||   (\bu_\la)^{p-1 }(\bu_\la- \bu) ||_l ||  (\bu_\la- \bu)^{q-1}   ||_l
\]
where $l = \frac{2Q}{Q+\al}$. From H\"{o}lder's inequality,
\[ \int\limits_{E_\la} (\bu_\la)^{l(p -1) }(\bu_\la- \bu)^l  \leq  (\int\limits_{E_\la} (\bu_\la)^{ls(p -1) } )^{1/s}  (\int\limits_{E_\la} (\bu_\la- \bu)^{s'l})^{1/s'}, \quad \frac 1s + \frac{ 1}{s'} = 1. \] 
Choosing $s$ such that $ ls(p -1) <  p+1$ and   $q  > \max\{ s'l, (q-1)l\}$   ( respectively, when $ p = \si$ choose  $s = \frac{Q+\al}{2\al}$ so that $s' = \frac{Q+\al}{Q-\al}$, and $q = \frac{l(\si+1) }{(\si+1) - l(\si-1)}$   so that $(q-1)l = q$ ), we get
\[||   (\bu_\la)^{p -1 }(\bu_\la- \bu) ||_l \leq || ( \bu_\la)^{p -1 } ||_{sl}||(\bu_\la- \bu) ||_{s'l} = ||  (\bu_\la)^{p -1 } ||_{sl} ||(\bu_\la- \bu) ||_{q}.
\]
After simplification,   
\bea\label{t2} \int\limits_{E_\la} \left(\bu_\la(\zeta) - \bu(\zeta)\right)^{q} \,d\zeta & \leq &  C \si ||   (\bu_\la)^{p -1 }||_{sl} ||(\bu_\la- \bu) ||_q ||  (\bu_\la- \bu)^{q-1}   ||_l \nonumber\\
& =&  
C \si  \left(  \int\limits_{E_\la}  (\bu_\la)^{(p +1) } \right)^{\al/Q} \left(\int\limits_{E_\la} \left(\bu_\la(\zeta) - \bu(\zeta)\right)^{q} \,d\zeta \right).
\eea
Since $\bu \in L^{p+1}(\HH^n)$, 
 \[\int\limits_{E_\la}  \bu_\la(\xi)^{(p +1) } d\xi = \int\limits_{{\mathcal R}^\la{E_\la}}  \bu(\xi)^{(p +1) } d \xi  \leq \int\limits_{{\Si_\la}^c}  \bu(\xi)^{(p +1) } d \xi < \frac 12 \mbox{ (say) for all } \la < -M_0
\] for some $M_0 >0$, where \be {\mathcal R}^\la{E_\la} =\{ (
\bar z, 2\la  -t) : (z, t) \in E_\la \}\ee is the $\HH$-reflection of the set $E_\la$ with respect to the plane $t = \la$ and $\Si_\la^c$ is complement of $\Si_\la$. 
For all such $\la < -M_0$, (\ref{t2}) will imply that 
$  \int\limits_{E_\la} \left(\bu_\la(\zeta) - \bu(\zeta)\right)^{q} \,d\zeta  = 0$ and hence $meas(E_\la)  = 0 $ for all $ \la < -M_0$. In particular,
from (\ref{sets}) we further conclude that
\be \bu_\la(\zeta\eta) \leq \bu (\zeta\eta ) \mbox{  and }  \bu_\la(\zeta\eta_\la) \leq \bu (\zeta\eta_\la ) \mbox{ for all  } \eta \in \Si_\la.\ee

If $u(\zeta) =  u(\zeta_\la)$ for some $\zeta \in \Si_\la^*$, then we get our symmetry (\ref{ksym}) as discussed in (iii). 

Otherwise, $u(\zeta) >  u(\zeta_\la)$ in $\Si_\la^*$ for all $\la < - M_0$.


 and  hence 
\be E_\la = \emptyset, \mbox{ the  empty set for all }  \la < - M_0 \ee
In this case,   define $ \La := \sup \{ \la < 0 :  \bu_\la(\zeta \xi)<\bu(\zeta \xi ) \mbox{ or }    \bu_\la(\zeta\xi_\la)<\bu (\zeta\xi_\la ) \, \mbox{ for all }\xi \in \Si_\la \mbox { and } \mbox{  for all } \zeta \in \Si_\la^*   \}$. 

{\bf Step 2: 	\bea &&  \bu_\La(\zeta)  \equiv   \bu(\zeta  )  \mbox{ for  all }  \zeta \in   \Si_\La^* \\ 
\mbox{OR }&&	 \bu_\La(\zeta \xi) \equiv  \bu(\zeta \xi )   \mbox { and } \bu(\zeta \xi_\La) \equiv \bu_\La(\zeta \xi_\La )\mbox{ for  some }  \zeta \in   \Si_\La^*  \mbox{ for all  } \xi \in \Si_\La.  \nonumber\\  \label{c1}\eea  }
Continuity of the map  $ \la \mapsto \bu(\zeta \xi) - \bu_\la(\zeta \xi ) $ in the $\la$ variable  implies  that
\be \label{c0}  \bu(\zeta \xi) \geq \bu_\La(\zeta \xi )  \,  \mbox { and } \bu(\zeta \xi_\la) \geq \bu_\La(\zeta \xi_\la )\mbox{  for all } \xi \in \Si_\La, \,\zeta \in \Si_\la^*.\ee
If for some $\zeta \in \Si_\La^*$, if  $\bu( \zeta) =  \bu_\La( \zeta)$ then we get the symmetry (\ref{c1}).
Otherwise, we must have
\be\label{la-contra} \bu( \zeta) >  \bu_\La( \zeta)  \mbox{  for all } \zeta \in \Si_\La^*.\ee
We will show that  (\ref{la-contra}) is not possible if $\La < 0$. 

Suppose $\La<0$ and (\ref{la-contra}) holds. If we denote $E_\La = \{ \xi \in \Si_\La^* : \bu_\La(\xi) >  \bu(\xi)\}$ and $\overline{E_\La}$ as its closure, then
$E_\La = \emptyset$ and  $meas(\overline{E_\La}) = 0$. Since we can write  
  $ \overline{E_\La} = \limsup\limits_{\la \to \La, \la > \La } E_\la$,  for any given $\var > 0$, there exists $\delta > 0$ such that  $\La + \delta < 0$ and for all $\la \in [\La, \La + \delta]$, 
\be     
  \int\limits_{E_\la}  \bu_\la(\xi)^{(\si +1) } d\xi < \var .\ee
 Repeating the arguments of Step 1, we conclude that $meas(E_\la) = 0 $ for all $\la \in [\La, \La + \delta]$ and hence, 
$E_\la = \emptyset $ for all $ \la \in [\La, \La + \delta]$ (as $E_\la$ is an open set). This   contradicts the definition of $\La$. 
  Hence, we must have either (\ref{c1}) holds for some $\zeta$ or that  $\bu \equiv \bu_\La$ in $\Si_\La^*$.

  If  $\La = 0$ then (\ref{c0}) implies 
   in particular that 
  \be \bu(\zeta ) \geq \bu_\la(\zeta  ), \,  \mbox{  for all } \zeta \in \Si_\la^* \mbox{ and } \mbox{  for all } \la \leq 0 . \label{inq} \ee
  Now, moving the plane from $\la >>0$ large, repeat the Steps 1 and 2 for $\bu -\bu_\la$ with $\la \geq 0$.  Due to (\ref{inq}) the process cannot stop for $\la > 0$  and we conclude
   \[\bu \leq \bu_\la \mbox{ in } \Si_\la^* \mbox{ for all }\la \geq 0.\] It follows that  $\bu$ is invariant with respect to the $\HH$-reflection about the plane $t= 0$,
\be \label{31}\bu(z, t) = \bu(\bar z, -t) \quad \mbox{ for all } \quad (z , t) \in \HH^n. \ee

{\bf Remark: \/} If $\La < 0$ then the symmetry with respect to plane $t= \La$ implies that $ \bu(0) = \bu_\la(0)$, so that singularity can be removed. But if $\La = 0$,   we cannot conclude any  more information about the singular behaviour of $\bu$ at the origin.

\bigskip
 {\bf Step 3:  $\HH$-symmetry \/}  \\
 Case (i): $\La = 0$  and  suppose that $\bu$ is singular at the origin: For $\La =0$, from  Step 2 above we have  
 \be \bu( z, t) = \bu(\bar z, -t) \quad \mbox{ for all } (z, t) \in \Si_0, \mbox{ i.e., } t \geq 0. \ee
 Repeating the process with the $\HH$-reflection
 \[ {\mathcal R}_x: (z, t) \mapsto (-\bar z, 2\la_2 -t)\]
 for the function $\bu$, we will find a point $\zeta_0$ and $\la_1 \leq 0 $  such that 
 \be \bu(\zeta_0 \xi)  = \bu_\la(\zeta_0 \xi)   \quad \mbox{ for all } \xi \in \Si_{\la_1}.\ee 
 Now since we have assumed that $\bu$ is singular at the origin, the point $\zeta_0$ is necessarily the origin and hence
 \be  \bu(z, t) = \bu(-\bar z, -t) =  \bu(- z, t) \quad \mbox{ for all } t \geq 0 .\ee
 In fact, since in this case we know that always $\La =  0$, and any direction can be chosen as $x$-direction after a rotation, it follows that 
 \be  \bu(z, t) = \bu( |z |, t)  \quad \mbox{ for all } t \geq 0 .\ee

 Case (ii): In general,  suppose that there exists $\zeta_1 \in  \Si_\la^*$  and $\la_1$ such that 
 \be \label{-0}\bu_{\la_1}(\zeta_1\eta) = \bu (\zeta_1\eta ) \mbox{  and }  \bu_{\la_1}(\zeta_1\eta_{\la_1}) = \bu (\zeta_1\eta_{\la_1} ) \mbox{ for all  } \eta \in \Si_{\la_1} \ee 
 where the $\HH$ reflection here is \[{\mathcal R}_y(z,t) \mapsto (\bar z, 2\la-t).\] 
 If $\zeta_1 = (z_1, t_1)$, $\xi = (z, t)$  then (\ref{-0}) implies that 
\bea 
\bu( z+ z_1, t+t_1+ 2 Im(z_1 \bar z)) &= &\bu( \overline{z+ z_1}, 2 \la_1 - (t+t_1+ 2 Im(z_1 \bar z))) \nonumber\\
&= &\bu( \overline{z+ z_1}, 2 \la_1 - t - (t_1+ 2 Im(z_1 \bar z)))\\
\bu( \bar{z} + z_1, 2\la_1 - t + t_1+ 2 Im(z_1 { z})) &= &\bu( \overline{\bar {z}+ z_1}, 2 \la_1 - ( 2\la_1 - t+t_0+ 2 Im(z_1  z))) \nonumber\\
&= &\bu(  {z}+ \bar{ z_1},  t - (t_1+ 2 Im(z_1  z)))
.\eea
As a second step towards proving symmetry of $\bu$,  we apply above steps with the $\HH$ reflection 
  \be {\mathcal R}_\theta : (z, t) \mapsto (e^{i\theta} \bar z, 2\la - t), \,\, \theta \in S^1, \,\, \theta \neq 0 \,\,e^{i\theta} \bar z =( e^{i\theta} \bar z_1, \ldots, e^{i\theta} \bar z_n) \in \CC^n. \ee 
Then, we will obtain a point $\zeta_2$ say, depending on $\theta$ and $\la_2$ such that
\be \bu_{\la_2}(\zeta_2\eta) = \bu (\zeta_2\eta ) \mbox{  and }  \bu_{\la_2}(\zeta_2\eta_{\la_2}) = \bu (\zeta_2\eta_{\la_2} ) \mbox{ for all  } \eta \in \Si_{\la_2}. \ee 
We claim that $ \zeta_2 = \zeta_1$ and $\la_1 = \la_2$:\\
Given a point $\eta = (z, t) \in \HH^n$, either $ \eta \in \zeta_1 \Si_{\la_1}$ or  $ \eta \in \zeta_1 \HH^n \setminus \Si_{\la_1}$. Hence, writing $\eta = \zeta_1 \xi$ if $ \eta \in  \zeta_1\Si_{\la_1}$  ( respectively, $\eta = \zeta_1 \xi_{\la_1}$ if $\eta \in \zeta_1 \HH^n \setminus \Si_{\la_1}$) then from (\ref{-0}) we conclude that 
\be\label{-21} \bu(\eta) = \bu(\eta_{\la_1}) \quad \mbox{ for all } \eta \in \HH^n.\ee
Similarly, splitting the space $\HH^n$ with respect to the hyper surface $\zeta_2{\mathcal H}_{\la_2}$ we conclude
\be\label{-22}  \bu(\eta) = \bu(\eta_{\la_2}) \quad \mbox{ for all } \eta \in \HH^n \ee
i.e., $\bu$ must be symmetric with respect to both the planes $t = \la_1$ and $ t= \la_2$ which is possible only if $\bu$ is constant, and hence $\bu = 0$ which is a contradiction since we assumed that $ \bu > 0$. Hence the claim follows.

Thus, if $ \zeta_1 = (z_1, t_1)$ with $z_1 \neq 0$ and $ \la_1$ is such that (\ref{-0}) holds  then necessarily $\zeta_1$ continues to remain the point of symmetry of $\bu$  for all the reflections ${\mathcal R}_\theta$, $\theta \in S^1$ 
\[{\mathcal R}_\theta^{\la_1}(z,t) \mapsto (\bar z, 2\la_1-t).\]and we conclude that
\be\bu(\zeta_1\eta) = \bu ( {\mathcal R}_\theta^{\la_1}    ( \zeta_1\eta) ) \mbox{  and }  \bu(\zeta_1   {\mathcal R}_\theta^{\la_1}    \eta  ) = \bu (   {\mathcal R}_\theta^{\la_1}    (  \zeta_1 {\mathcal R}_\theta ^{\la_1}   \eta )) \mbox{ for all  } \eta \in \Si_{\la_1} \ee 
This implies that
\be  \label{zeta-cyl0} \bu(\zeta_1\eta) = \bu ( {\mathcal R}_\theta^{\la_1}    ( \zeta_1\eta) ) \mbox{ for all  } \theta \in S^1, \ee
i.e., 
\be\label{zeta-cyl} \bu(z+ z_1, t+t_1+ 2 Im(z_1 \bar z))      = \bu (e^{i\theta}(z+ z_1), 2 \la_1-( t+t_1+ 2 Im(z_1 \bar z)) ) \mbox{ for all  } \theta \in S^1. \ee
If $z_1 = 0$ then
\be \bu(z, t) = \bu(e^{i\theta} z, 2\la_1 -t)  \mbox{ for all  } \theta \in S^1 \ee implies that $\bu$ is cylindrical and symmetric in $t$-variable with respect to the plane $t = \la_1$.

\qed

{\bf Step 4:  The limit $ \lim\limits_{|\xi| \to \infty} |\xi|^{Q-\al} u (\xi)$ exists for $p = \si$ } \\
If in Step 2,  $ \La < 0$ then we can define $\bu(0) = \bu(0_\La) = \bu( 0, 0, \La)$ and hence $\bu$ has no singularity at origin, i.e.,  the limit  \be \lim\limits_{|\xi| \to \infty} |\xi|^{Q-\al} u (\xi)  = u_\infty \quad \mbox{  exists} \ee
as \be \lim\limits_{|\xi| \to \infty} |\xi|^{Q-\al} u (\xi)  = \lim\limits_{|\xi| \to 0} \bu (\xi)
=\lim\limits_{|\xi_\La| \to 0_\La} \bu (\xi_\La) =\bu(0_\La).\ee
For $\La = 0$, $u$ is $\HH$-symmetric with respect to the plane $t = 0$.  
If $\lim\limits_{|\xi| \to \infty} |\xi|^{Q-\al} u (\xi)$ is not finite, then  perform the  CR transform of $u$  with respect to any point $(0,0, t_0)$ on the $t$ axis i.e., define
\be \bu_{t_0}( z, t) :=  \bu( z, t+t_0) \ee
so that $\bu_{t_0}$ has singularity at  $(0,0, t_0)$. Then repeating Step 1 and 2 for $\bu_{t_0}$ we conclude that $\bu_{t_0}$ is $\HH$-symmetric with respect to the plane $t = t_0$. Hence $\bu$ is symmetric with respect to the plane $t = t_0$. Since $t_0$ was arbitrary, $\bu$ is $\HH$-symmetric with respect to the plane $t=t_0$ for all $t_0 \in \R$. This is possible only if $\bu$ is independent of the $t$ variable and hence $u$ must also be independent of the $t$ variable. 

Since $u$ satisfies (\ref{m1}), it follows that $u \equiv 0$. For if $\xi_1 = (0, 0)$, $\xi_2 = (0, - t_1)$ with $ 0 < t_1 $  are points on the $t$ axis, then $ u(\xi_1) = u(\xi_2)$ implies that 
\be  \int\limits_{\HH^n} G(\xi_1, \eta) u(\eta)^p \, d\eta  = \int\limits_{\HH^n} G(\xi_2, \eta) u(\eta)^p \, d\eta\,\,\, i.e., \int\limits_{\HH^n} [G(\xi_1, \eta) - G(\xi_2, \eta)] u(\eta)^p \, d\eta = 0 \ee a contradiction since $G(\xi_1, \eta) = |(z,t)|^{-(Q-\al)} > |(z, t +t_1)|^{-(Q-\al)} =  G(\xi_2, \eta)$ for any $\eta =(z,t) \in \HH^n$. 
Thus either $u \equiv 0$ or if  $u > 0$ then the limit 
\[ \lim\limits_{|\xi| \to \infty} |\xi|^{Q-\al} u (\xi) = u_\infty \mbox{ exists }.\]
\qed

\subsection{Symmetry of solutions of (\ref{km3}) for the subcritical case $p < \frac{Q+\al}{Q-\al}$}\label{subcriti} As the proof here is similar to the critical case $p = \si$,  we only list here the main computations and discussions which allow us to arrive at the conclusion of cylindrical symmetry of solutions of (\ref{km3}).   Using the representation (\ref{i2}) we write 
 \bea\label{ii2}
&& \bu(\zeta_\la ) - \bu (\zeta) \nonumber\\
 & = & \int\limits_{\Si_\la } |\eta|^{-(Q-\al)} \left(f_\la(   \zeta\eta  ) - f(\zeta \eta) \right)d\eta 
 +  \int\limits_{\Si_\la }  | \eta_\la |^{-(Q-\al)} \left( f_\la(\zeta\eta_\la) -  f(\zeta \eta_\la) \right)  d\eta \nonumber\\
 & \leq & \int\limits_{A_{\la, \zeta} } |\eta|^{-(Q-\al)} \left(f_\la(   \zeta\eta  ) - f(\zeta \eta) \right)d\eta 
 +  \int\limits_{B_{\la, \zeta} }  | \eta_\la |^{-(Q-\al)} \left( f_\la(\zeta\eta_\la) -  f(\zeta \eta_\la) \right)  d\eta
 \eea 
 where as before \be
 E_\la:= \{ \xi \in \Si_\la^*: \bu_\la(\xi) > \bu(\xi)\} \ee
 and for $\zeta \in \Si_\la^*$, 
 \be\label{seta}
 A_{\la, \zeta} := \{\xi \in \Si_\la : \zeta \xi \in E_\la \} \mbox{ and  } 
 B_{\la, \zeta} := \{\xi \in \Si_\la : \zeta \xi_\la \in E_\la  \}. \ee
The arguments of Step I and Step 2 goes through once we write  
 \Bea  f(\xi_\la) - f(\xi) & < &  \frac{1}{|\xi_\la |^{(Q+\al) - p(Q-\al)}}(\bu_\la(\xi)^{p} - \bu(\xi)^{p}) \\
 &= &\frac{p a(\xi)^{p -1 } }{|\xi_\la |^{(Q+\al) - p(Q-\al)}} (\bu_\la(\xi)- \bu(\xi))\\
 & \leq &p \frac{\bu_\la(\xi)^{p -1 }}{|\xi_\la |^{(Q+\al) - p(Q-\al)}}  (\bu_\la(\xi)- \bu(\xi)). \Eea Here $a(\xi)$ is a real number between $\bu(\xi_\la)$ and $\bu(\xi)$ and  we use that  $a(\xi) < \bu_\la(\xi)$ for $\xi \in E_\la$. 
From the Claim in Step 1,  $meas(E_\la) = 0$ in sub critical case as well and from Step 2, there exists $\La \leq 0$ such that 
	\bea &&  \bu_\La(\zeta)  \equiv   \bu(\zeta  )  \mbox{ for  all }  \zeta \in   \Si_\La^* \\ 
\mbox{OR }&&	 \bu_\La(\zeta \xi) \equiv  \bu(\zeta \xi )   \mbox { and } \bu(\zeta \xi_\La) \equiv \bu_\La(\zeta \xi_\La )\mbox{ for  some }  \zeta \in   \Si_\La^*  \mbox{ for all  } \xi \in \Si_\La.  \eea  }
 
 In the Step 3,  $\La < 0$ is not possible. For if   $\la < 0$ and  $ \bu(\zeta) = \bu(\zeta_\la)$  for some $\zeta \in \Si_\la^*$ then  $f_\la(\zeta \xi ) =  f(\zeta \xi )$ for all $\xi \in \Si_\la$.  Since $\la < 0$, $0 \in \Si_\la$ and hence $f_\la(\zeta ) =  f(\zeta )$ i.e., 
\[ \frac{\bu(\zeta_\la)^{p} }{ |\zeta_\la|^{(Q+\al) - p(Q-\al)}}  = \frac{\bu(\zeta)^{p} }{ |\zeta|^{(Q+\al) - p(Q-\al)}} \mbox{ implying that  } |\zeta_\la|^{(Q+\al) - p(Q-\al)} = |\zeta|^{(Q+\al) - p(Q-\al)} \]
which is not true if $\la < 0$.  Hence, $\La = 0$ and we conclude (\ref{31}) as above. Case (i) of Step 3 now implies that $\bu$ is cylindrical and hence so is $u$.

\subsection{Nonexistence of solutions for $1 < p < \frac{Q+\al}{Q-\al}$ }\label{nonex}
{\bf Proof of Theorem \ref{non-exist}:\/}
Note that in case of sub-critical exponent i.e., $ 1 < p < \frac{Q+\al}{Q-\al} $,
 necessarily the singularity in the integral equation continues to persist for $\bu$ due to the term \[ \dfrac{1}{ |\xi|^{(Q+\al) - p(Q- \al)}   }\] forcing that $\La = 0$.  From the previous subsection, we  conclude that $u$ is $\HH$-symmetric with respect to the plane $t = 0$.  Now as 
explained  in the Step 4,  performing CR inversion with respect to any point $(0, t_0)$ on the $t$-axis and repeating above steps we get that $u$ must be cylindrical about the point $(0, t_0)$.    It follows that  $u$ must be  independent of the $t$ variable and  hence  $u \equiv 0$. 

\qed

\section{Classification of solutions of (\ref{mcrit}) }\label{unq}



Recall that we  refer to 
\be \label{std}
u_0 =  C_0 ( t^2 + (r^2 + 1)^2)^{- \frac{Q-\al}{4} } = C_0 |\om +i|^{\frac{Q-\al}{2}}\ee  as the  standard solution  of (\ref{mcrit}).
The following properties of $u_0$ can be verified easily:\\
(i) $ u_0(0) = C_0 = (u_0)_\infty = \lim\limits_{|\xi| \to \infty} |\xi|^{Q-\al} u_0(\xi) $, and \\
(ii) $u_0(\xi) = \bu_0(\xi)$ for all $\xi \in B(0,1)$, i.e., $u_0$ is inversion symmetric with respect to the unit CC-sphere $\pa B(0,1)$. \\
To classify the solutions of (\ref{mcrit}), we will first prove that (i) and (ii) are characteristic properties of the solutions of the integral equation (\ref{mcrit}), i.e.,  of
\be 
\label{m}
u (\zeta ) = \int\limits_{\HH^n} G_\al(\zeta, \xi)  u(\xi)^{\frac{Q+\al}{Q-\al}}\, d\xi ,\quad  0< \al <Q \ee
where \be  G_\al(\zeta, \xi) = 
|\zeta^{-1} \xi|^{-(Q-\al)} . \ee 
For $u$ a cylindrical solution of (\ref{m}), we define \be\label{sinv}
	 \bu_{s^2} (\xi) :=   \frac{ s^{Q-\al}}{\rho^{Q-\al}} u \left(\frac{s^2 r}{\rho^{2}}, - \frac{ s^4 t}{\rho^{4}}\right) = s^{Q-\al}(u \circ \delta_{s^2})^\Delta (\xi), ~~~s> 0 \ee which is the CR type inversion with respect to the CC sphere $\pa B(0, s)$. Here $(u \circ \delta_{s^2})^\Delta $ denotes the CR-type inversion of the function $u \circ \delta_{s^2}$ and $\delta_{s^2}$ is the dilation  $(z,t) \mapsto (s^2 z, s^4 t)$. The following lemma relates $u$ and $\bu_{s^2}$.


\begin{lemma}\label{iinv}
	The map $\bu_{s^2}$ is  a solution of (\ref{m}) and the integral (\ref{m}) can be expressed as 
	 \be\label{uvs}
	 u_s(\zeta) = \int\limits_{B(0,1)} G_\al(\zeta, \xi)  u_s(\xi)^{\si} \,d\xi  + \int\limits_{B(0,1)} G_\al( \hat\zeta, \xi) \frac{1}{ |\zeta|^{Q-\al}} \bu_{s}(\xi)^{\si} \, d\xi. \ee
Similarly,
\be\label{vsu} 
 \bu_{s}(\zeta) = \int\limits_{B(0,1)} G_\al(\zeta, \xi)  \bu_{s}(\xi)^{\si} \,d\xi  + \int\limits_{B(0,1)} G_\al( \hat \zeta, \xi) \frac{1}{ |\zeta|^{Q-\al}} u(\xi)^{\si} \, d\xi. \ee
	\end{lemma}
\proof  From Lemma \ref{ing} since $u_s$ is also a solution of (\ref{m}), it suffices to prove (\ref{uvs}) for $s= 1$, i.e., we prove that and solution of (\ref{m}) satisfies
\be\label{uv1}
u(\zeta) = \int\limits_{B(0,1)} G_\al(\zeta, \xi)  u(\xi)^{\si} \,d\xi  + \int\limits_{B(0,1)} G_\al( \hat\zeta, \xi) \frac{1}{ |\zeta|^{Q-\al}} \bu(\xi)^{\si} \, d\xi. \ee
Since $u$ solves (\ref{m}), we can write 
\[ u(\zeta) = \int\limits_{B(0,1)} G_\al(\zeta, \xi)  u(\xi)^{\si}\,d\xi  + \int\limits_{\HH^n \setminus B(0,1)} G_\al(\zeta, \xi)  u(\xi)^{\si}\,d\xi = I_1 +I_2\,  (\mbox{ say}). \]

In $I_2$, make the change of variables $ \eta = \hat \xi$ so that $| \eta| \leq 1$,  $d \xi = \dfrac{1}{|\hat \xi|^{2Q}} d \hat\xi = \dfrac{1}{|\eta|^{2Q}} d \eta $. Also, using  cylindrical symmetry we get
 $ u( \xi)  = u(z,t) = u(-z, t)  = u(\hat{\hat\xi}  ) = | \eta|^{Q-\al} \bu  ( \eta)$ we have 
\Bea
I_2 & = & \int\limits_{|\xi|> 1 } G_\al(\zeta, \xi)  u(\xi)^{\si}\,d\xi
 =  \int\limits_{|\eta|< 1 } G_\al(\zeta, \eta)  \{ | \eta|^{Q-\al} \bu  ( -\eta)\}^{\si} \, \dfrac{1}{|\eta|^{2Q}} d \eta \\
&  = & \int\limits_{|\eta|< 1 } G_\al(\zeta, \hat \eta)\dfrac{1}{|\eta|^{Q-\al}} \bu  ( \eta)^{\si} \, d\eta 
= \int\limits_{|\eta|< 1 } G_\al(\hat\zeta,  \eta)\dfrac{1}{ |\zeta|^{Q-\al}} \bu  ( \eta)^{\si} \, d\eta
\Eea
and (\ref{uv1}) follows. 

\qed

Next, we show  that any cylindrical  solution of (\ref{m}) is  inversion symmetric with respect to a CC-sphere $\pa B(0, s_0)$ for some $0< s_0 \leq 1 $.
\begin{Pro}\label{inv-sym}
	If $u$ is a cylindrical solution of (\ref{m}) and $s > 0$ be such that 
	\be s^{Q-\al} = \dfrac{u_\infty}{u(0)}. \ee
Then, 
	\be \label{inv}
	u \equiv   \bu_{s^2} \mbox { in } \HH^n.\ee 
In particular,  $ u(0) = u_\infty$ iff $s_0 = 1$ and $u \equiv \bu $ in $\HH^n$. 
\end{Pro}
\proof 
Let $s = 1$ so that $u_\infty = u(0)$ and $u$ is cylindrical about origin as well as even in the $t$ variable. As in \cite{CLO}, define   
\be v(z, t) =   (u \circ \tau_{(0,1)})^\Delta(z, t)=    \frac{1}{\rho^{Q-\al}} u(\frac{z}{\om}  , \frac{ -t  }{\rho^4} +1 ),   \, \rho^4 = | z|^4 + t^2.\ee
where $\bu$ is cylindrical, even in $t$ and $\bu(0)=u_\infty= u(0) = \bu_\infty$.
 Now $v$ is also a solution of (\ref{m}) as $\bu$  and a translation of $\bu$ is again a solution of (\ref{m}). Hence from our Theorem \ref{symm}, $v$ must be radial in $z$ variable with respect to some point $(z_0, t_0)$  in $\HH^n$.  
  It can be verified that 
$v(0,1) = u(0) = u_\infty$ and $v(0) = u_\infty= u(0) $.  Therefore $v$ must be $\HH$-symmetric  about the point $ O = (0, 1/2) \in \CC^n \times \RR$.
   Therefore, $v$ satisfies (\ref{zeta-cyl}) with $ \zeta = (0, 1/2)$ and we get  
   \be 
   v(z, t) = v(e^{i\theta} \bar{z}, 1-t) \quad \mbox{ for all } (z, t) \in \HH^n. \ee 
   Therefore,
\bea
\frac{1}{\rho^{Q-\al}} u (\frac{z}{\om}  , \frac{ -t  }{\rho^4} -1 ) & = &  \frac{1}{|(z, 1-t)|^{Q-\al}} u (\frac{e^{i \theta} \bar{z}}{\tilde\om}  , \frac{-1 +t  }{|\tilde \om|^2}  +1 ), \mbox{ where } \tilde\om = (1-\bar \om),  \om = t + i |z|^2, \nonumber\\
  \mbox{ i.e., }  
u (\frac{z}{\om}  , \frac{ -t  }{|\om|^2} -1 ) & = &  \frac{|\om|^{(Q-\al)/2}}{|1-\bar{\om}|^{(Q-\al)/2}} u (\frac{e^{i \theta} \bar{z}}{(1-\bar \om)}  , \frac{-1 +t  }{|1-\bar \om|^2}  +1 )\nonumber \\
&=& \frac{|\om|^{(Q-\al)/2}}{|1-\bar{\om}|^{(Q-\al)/2}} u (\frac{ |z|}{|1-\bar \om|}  , \frac{-1 +t  }{|1-\bar \om|^2}  +1 ) \label{inv00}
\eea
   since $u$ is cylindrical. 
   Introduce variables $z_1 = \frac{|z|}{|\om|}=   \frac{|z|}{\rho^2}$, 
$t_1 = \frac{ -t  }{\rho^4} -1  = \frac{ -t  }{|\om|^2} -1 $, we have 
\[ |\om_1 = t_1 + i |z_1|^2| = \frac{|1 - \bar \om|}{|\om|}\] and we see from (\ref{inv00}) that 
\be
u(|z_1|, t_1) = \frac{1}{|\om_1|^{(Q-\al)/2}} u ( \frac{ |z_1|}{| \om_1|}  , \frac{t_1}{|\om_1|^2} )
\ee 
   which is the required result. 

\qed

{\bf Remark:\/} Due to the cylindrical  symmetry  proved in the previous section, the conclusion of   Proposition \ref{inv-sym} also follows from the results of \cite{monti}.

\begin{corollary}\label{3.1}
	If $u$ is a cylindrical solution of (\ref{m1}) with  \be
	\label{s} s^{Q-\al} = \frac{u_\infty}{u(0)}\ee then 
	\be \label{l3.1}
	u(sr, s^2t) = \frac{1}{\rho^{Q-\al}} u\left(\frac{s r}{\rho^{2}}, \frac{s^2 t}{\rho^{4}}\right). \ee
	In general, for any $\xi_0 \in \H$, defining $s^{Q-\al} = \dfrac{u_\infty}{u(\xi_0)}$ we have 
	\be \label{l3.1a}
	u\circ \tau_{\xi_0}(\delta_s \xi ) = \frac{1}{\rho^{Q-2}} u\circ \tau_{\xi_0}( \delta_s \hat \xi ). \ee
	
\end{corollary}
\proof Define $u_1(z, t) = s^{\frac{Q-\al}{2}}u(sz, s^2t)$. Then $u_1(0) = u_\infty$ and $\lim\limits_{\rho \to \infty} \rho^{Q-2}u_1(z,t) = u_\infty$. Thus $u_1$ satisfies the conditions of Lemma \ref{inv-sym} and hence 
\[
u_1(r,t) = \frac{1}{\rho^{Q-2}} u_1(\frac{r}{\rho^{2}}, \frac{t}{\rho^{4}}) \]
which proves (\ref{l3.1}).  Arguing similarly for the translated function $v( \xi) = u\circ \tau_{\xi_0} (\xi) = u( \xi_0 \cdot \xi)$, we get (\ref{l3.1a}). 

\qed

\begin{lemma}\label{uniq}
If $u$, $w$ are two cylindrical solutions of the integral equation (\ref{m}) $\HH$-symmetric with respect to the hyperplane through the origin such that 
\be\label{10}  u_\infty = w_\infty \ee
then $ u \equiv w$.
\end{lemma}

\proof 


Without loss of generality, let  $u(0) = u_\infty = w_\infty$ and suppose that 
\be \label{contra}
 w(0) > u(0). \ee 
For $\la > 0$,  let
\Bea
\bu_{\la^2} (\xi) & := &  \frac{ \la^{Q-\al}}{\rho^{Q-\al}} u \left(\frac{\la^2 r}{\rho^{2}}, - \frac{ \la^4 t}{\rho^{4}}\right)  \\
\bw_{\la^2} (\xi) & := &  \frac{ \la^{Q-\al}}{\rho^{Q-\al}} w \left(\frac{\la^2 r}{\rho^{2}}, - \frac{ \la^4 t}{\rho^{4}}\right)  \Eea so that 
\bea
\bu_{\la^2}(\zeta) = \int\limits_{B(0,\la)} G_\al(\zeta, \xi)  \bu_{\la^2}(\xi)^{\si} \,d\xi  + \int\limits_{B(0,\la)} G_\al( - \delta_{\la^2}\hat \zeta, \xi) \frac{\la^{Q-\al}}{ |\zeta|^{Q-\al}} u(\xi)^{\si} \, d\xi \\
\mbox{ and }  
\bw_{\la^2}(\zeta) = \int\limits_{B(0,\la)} G_\al(\zeta, \xi)  \bw_{\la^2}(\xi)^{\si} \,d\xi  + \int\limits_{B(0,\la)} G_\al( - \delta_{\la^2}\hat \zeta, \xi) \frac{\la^{Q-\al}}{ |\zeta|^{Q-\al}} w(\xi)^{\si} \, d\xi. \eea
Then,
\bea \label{uw}
\bu_{\la^2}(\zeta) - \bw_{\la^2}(\zeta)& =& \int\limits_{B(0,\la)} G_\al(\zeta, \xi)  (\bu_{\la^2}(\xi)^{\si} - \bw_{\la^2}(\xi)^\si)\,d\xi  \nonumber\\
&& \quad+ \int\limits_{B(0,\la)} G_\al( - \delta_{\la^2}\hat \zeta, \xi) \frac{\la^{Q-\al}}{ |\zeta|^{Q-\al}} (u(\xi)^{\si}-w(\xi)^\si ) \, d\xi. \eea
Define 
\be 
A_\la := \{ \xi \in \overline{ B(0, \la)} : \bu_{\la^2}(\xi) > \bw_{\la^2}(\xi) \}. \ee 
Since $w(0) > u(0)$ there exists $\la_0$ such that for all $|\xi| < \la_0$, \[ w(\xi) > u(\xi).\]
Hence from (\ref{uw}), for all $\la \leq \la_0$, 
\be
\bu_{\la^2}(\zeta)  - \bw_{\la^2} (\zeta) \leq  \int\limits_{B(0,\la)} G_\al(\zeta, \xi) ( \bu_{\la^2}(\xi)^{\si}  - \bw_{\la^2}(\xi)^{\si} ) \,d\xi \leq \int\limits_{A_\la} G_\al(\zeta, \xi) ( \bu_{\la^2}(\xi)^{\si}  - \bw_{\la^2}(\xi)^{\si} ) \,d\xi. \ee 
For $\zeta \in A_\la$, multiplying the above equation by $(\bu_{\la^2}(\zeta)  - \bw_{\la^2}(\zeta))^{q-1}$   where $q = \frac{2Q}{Q-\al}$ and applying HLS inequality together with H\'{o}lder's inequality as in Step 1 of the proof of Theorem \ref{symm}  we get
\bea\label{la1}
&&\int\limits_{A_\la} ( \bu_{\la^2}(\zeta) - {\bw_{\la^2}}(\zeta))^{q} \, d\zeta \nonumber\\
& \leq & \si \int\limits_{A_\la}\int\limits_{A_\la} G_\al(\zeta, \xi) \bu_{\la^2}(\xi)^{\si-1} (\bu_{\la^2}(\xi)  - \bw_{\la^2}(\xi) ) (\bu_{\la^2}(\xi)  - \bw_{\la^2}(\xi) )^{q-1} \, d\zeta \,d\xi \nonumber \\
& \leq & C\si || \bu_{\la^2}(\xi)^{\si-1} (\bu_{\la^2}(\xi)  - \bw_{\la^2}(\xi) )||_{L^l({A_\la})} || (\bu_{\la^2}(\xi)  - \bw_{\la^2}(\xi) )^{q-1} ||_{L^l({A_\la})} 
\mbox{  with } l = \frac{2Q}{Q+\al} \nonumber \\
& \leq &   C \si||u||_\infty^{\al/Q} meas(A_\la) \int\limits_{A_\la} ( \bu_{\la^2}(\zeta) - {\bw_{\la^2}}(\zeta))^{q} \, d\zeta\eea
where choosing  $s = \frac{Q+\al}{2\al} $ so that $s' = \frac{Q+\al}{Q-\al}$ and $s'l = \frac{2Q}{Q-\al} = q$
\Bea  || \bu_{\la^2}(\xi)^{\si-1} (\bu_{\la^2}(\xi)  - \bw_{\la^2}(\xi) )||_{L^l{A_\la})} &\leq & || \bu_{\la^2}(\xi)^{\si-1}||_{L^{sl}(A_\la)} || (\bu_{\la^2}(\xi)  - \bw_{\la^2}(\xi) )   ||_{L^{s'l} (A_\la)} \\
&\leq & C (||u||_\infty)^{\al/Q}) meas(A_\la)  || \bu_{\la^2}(\xi)  - \bw_{\la^2}(\xi)   ||_{L^{q} (A_\la)}.   \Eea
For $\la < \la_0$ small, $meas(A_\la)$ can be chosen small enough so that we get a contradiction in the inequality (\ref{la1}). Hence  
 we can  conclude that  \bea && meas(A_\la)  =  0 \mbox{ for all  } \la \leq  \la_0 \\ 
&& \mbox{ i.e., }  \bu_{\la^2}(\xi) \leq  \bw_{\la^2}(\xi) \mbox{ for } \xi \in B(0, \la) \mbox{ and for all sufficiently small } \la \leq  \la_0  \eea 
Now let $ \La = \sup\{\la \leq 1 : meas(A_\la) = 0 \}$. Then $meas(A_\La) = 0$.

We next observe that 
\be meas(B_\la) := \{ \xi \in B(0, \la) : u(\xi) > w(\xi) \} = 0 \mbox{ for all } \la \leq \La. \ee 
To see this we write 
\Bea u(\zeta) & = & \int\limits_{B(0,\la)} G_\al(\zeta, \xi)  u(\xi)^{\si} \,d\xi  + \int\limits_{B(0,\la)} G_\al( - \delta_{\la^2}\hat\zeta, \xi) \frac{\la^{Q-\al}}{ |\zeta|^{Q-\al}} \bu_{\la^2}(\xi)^{\si} \, d\xi\\
 w(\zeta) &= &\int\limits_{B(0,\la)} G_\al(\zeta, \xi)  w(\xi)^{\si} \,d\xi  + \int\limits_{B(0,\la)} G_\al( - \delta_{\la^2}\hat\zeta, \xi) \frac{\la^{Q-\al}}{ |\zeta|^{Q-\al}} \bw_{\la^2}(\xi)^{\si} \, d\xi\Eea and hence 
\bea 
 u(\zeta) - w(\zeta)  &= & \int\limits_{B(0,\la)} G_\al(\zeta, \xi) ( u(\xi)^{\si} - w(\xi)^\si )\,d\xi  \nonumber\\
&&\quad+ \int\limits_{B(0,\la)} G_\al( - \delta_{\la^2}\hat\zeta, \xi) \frac{\la^{Q-\al}}{ |\zeta|^{Q-\al}} ( \bu_{\la^2}(\xi)^{\si}  - \bw_{\la^2}(\xi)^\si )\, d\xi 
 \label{-3}\\
  & \leq & \int\limits_{B(0,\la)} G_\al(\zeta, \xi) ( u(\xi)^{\si} - w(\xi)^\si) \,d\xi ~~~\mbox{ for all } \la \leq \La. \eea
Again, arguing as in (\ref{la1}) for $u-w$, we get $meas(B_\la) = 0$ for all $\la \leq \La$. Since $\bu_{\la^2} - \bw_{\la^2}$ satisfies (\ref{uw})  and $u-w$ satisfies (\ref{-3}), we further conclude that
\bea && \mbox{ both } u < w \mbox{ and } \bu_{\la^2} < \bw_{\la^2} \mbox{ in } B(0, \la) \label{-4}\\
\mbox {or } && \mbox{ both } u \equiv w \mbox{ and } \bu_{\la^2} \equiv \bu_{\la^2} \mbox{ in } B(0, \la). \label{-5} \eea
Due to our assumption (\ref{contra}), the equality (\ref{-5}) is not possible and hence (\ref{-4}) must hold. 

  From Proposition \ref{inv-sym},  $w$ is inversion symmetric with respect to a CC-sphere $\pa B(0, s)$ of radius, say, $s_0$, i.e.,
\be  w(\xi) = \bw_{s_0^2}(\xi) = \frac{s_0^{Q-\al}}{\rho^{Q-\al}} w (s_0^2\hat \xi) \quad \mbox{ for all } \xi \in \HH^n. 
 \ee
 Since $s_0^{Q-\al} = \frac{w_\infty}{w(0)}$, \[w(0) = s_0^{Q-\al} w_\infty = s_0^{Q-\al} u(0) > u(0),\] we must have $s_0  > 1$. 
Also, $u$ is inversion symmetric with respect to the unit CC sphere.   We claim that $ \La = 1$. For if $ \La < 1$, then (\ref{-4}) holds so that both  $meas(A_\La) = 0$ and $meas(B_\la) = 0$. Given $\var > 0$ choose $ \delta > 0$ small with $\La+\delta < 1$ such that  $ meas(A_\la) < \var$ and $meas(B_\la) < \var$ for all $ \la \in [\La, \La + \delta]$. Now from (\ref{uw}) for $\la \in [\La, \La +\delta]$, 
\Bea
\bu_{\la^2}(\zeta) - \bw_{\la^2}(\zeta) &\leq & \int\limits_{A_\la} G_\al(\zeta, \xi)  (\bu_{\la^2}(\xi)^{\si} - \bw_{\la^2}(\xi)^\si )\,d\xi  \\
&&\quad+ \int\limits_{B_\la} G_\al( - \delta_{\la^2}\hat \zeta, \xi) \frac{\la^{Q-\al}}{ |\zeta|^{Q-\al}} (u(\xi)^{\si}-w(\xi)^\si ) \, d\xi 
\\
& \leq & \si \int\limits_{A_\la} G_\al(\zeta, \xi) \bu_{\la^2}(\xi)^{\si-1} (\bu_{\la^2}(\xi) - \bw_{\la^2}(\xi))\,d\xi  
\\
&&+ \si \int\limits_{B_\la} G_\al( - \delta_{\la^2}\hat \zeta, \xi) \frac{\la^{Q-\al}}{ |\zeta|^{Q-\al}}u(\xi )^{\si-1} (u(\xi)-w(\xi) )\, d\xi. \Eea
Multiplying above inequality by  $(\bu_{\la^2}(\zeta)  - \bw_{\la^2} (\zeta))^{q-1}$   where $q = \frac{2Q}{Q-\al}$,
\bea
&&\int\limits_{A_\la} |\bu_{\la^2}(\zeta) - \bw_{\la^2}(\zeta)|^q \,d\zeta \nonumber\\
  &\leq & \si \int\limits_{A_\la} \int\limits_{A_\la} G_\al(\zeta, \xi) \bu_{\la^2}(\xi)^{\si-1} (\bu_{\la^2}(\xi) - \bw_{\la^2}(\xi))(\bu_{\la^2}(\zeta)  - \bw_{\la^2} (\zeta))^{q-1}\,d\xi  \,d\zeta \label{int1}\\
&& + \si \int\limits_{A_\la}  \int\limits_{B_\la} G_\al( - \delta_{\la^2}\hat \zeta, \xi) \frac{\la^{Q-\al}}{ |\zeta|^{Q-\al}}u(\xi )^{\si-1} (u(\xi)-w(\xi) ) (\bu_{\la^2}(\zeta)  - \bw_{\la^2} (\zeta))^{q-1}\, d\xi \, d\zeta \label{int2}
\eea
and  the first integral on the RHS can be estimated as in (\ref{la1}) to get
\bea 
&&\si \int\limits_{A_\la} \int\limits_{A_\la} G_\al(\zeta, \xi) \bu_{\la^2}(\xi)^{\si-1} (\bu_{\la^2}(\xi) - \bw_{\la^2}(\xi))(\bu_{\la^2}(\zeta)  - \bw_{\la^2} (\zeta))^{q-1}\,d\xi  \,d\zeta \nonumber\\
&\leq & C \si||u||_\infty^{\al/Q} meas(A_\la) \int\limits_{A_\la} ( \bu_{\la^2}(\zeta) - {\bw_{\la^2}}(\zeta))^{q}. \eea
The second integral is 
\Bea
&&\si \int\limits_{A_\la}  \int\limits_{B_\la} G_\al( - \delta_{\la^2}\hat \zeta, \xi) \frac{\la^{Q-\al}}{ |\zeta|^{Q-\al}}u(\xi )^{\si-1} (u(\xi)-w(\xi) ) (\bu_{\la^2}(\zeta)  - \bw_{\la^2} (\zeta))^{q-1}\, d\xi \, d\zeta \\
&=& \si \int\limits_{A_\la}  \int\limits_{B_\la} G_\al( - \delta_{\la^2}\hat \zeta, \xi) \frac{\la^{Q-\al}}{ |\zeta|^{Q-\al}}u(\xi )^{\si-1} [u(\xi)-w(\xi) ] \frac{\la^{Q+\al}}{|\zeta|^{Q+\al}}[u(-\delta_{\la^2}\hat\zeta)  - w (-\delta_{\la^2}\hat\zeta)]^{q-1}\, d\xi \, d\zeta \\
 &=& \si \int\limits_{\widehat{A_\la}}  \int\limits_{B_\la} G_\al( - \delta_{\la^2}\hat \zeta, \xi) u(\xi )^{\si-1} [u(\xi)-w(\xi) ] [u(-\delta_{\la^2}\hat\zeta)  - w (-\delta_{\la^2}\hat\zeta)]^{q-1}\, d\xi \, d( \delta_{\la^2}\hat\zeta)  \\
 && \quad \mbox{ where } \widehat{A_\la} := \{ \delta_{\la^2}\hat \zeta \in \HH^n : \zeta  \in A_\la \}\\
&\leq & C \si ||   u^{\si-1 }( u - w) ||_{L^l(B_\la)} ||  (\bu_{\la^2}-  \bw_{\la^2})^{q-1}   ||_{L^l({A_\la})} 
\mbox{  with } l = \frac{2Q}{Q+\al}
\Eea
where we used the fact that 
\[\int\limits_{\widehat{A_\la}}  |[u(-\delta_{\la^2}\hat\zeta)  - w (-\delta_{\la^2}\hat\zeta)]^{q-1}|^l  d( \delta_{\la^2}\hat\zeta) =  \int\limits_{A_\la} |(\bu_{\la^2}-  \bw_{\la^2})^{q-1}|^l \, d\zeta.  \]

 From H\"{o}lder's inequality, choosing $s = \frac{Q+\al}{2\al} $ so that $s' = \frac{Q+\al}{Q-\al}$ and $s'l = \frac{2Q}{Q-\al} = q$,  we get 
\be\label{a1} ||   u^{\si-1 }( u - w) ||_{L^l(B_\la)} \leq ||u^{\si-1 }||_{L^{sl}(B_\la)} ||( u - w) ||_{L^{s'l}(B_\la)} = ||u^{\si-1 }||_{L^{sl}(B_\la)} ||( u - w) ||_{L^{q}(B_\la)} \ee
Moreover,  for $\la \in [\La, \La + \delta]$, 
 \be \label{a2} C \si  \left(\int\limits_{B_{\la}}  u^{(\si +1) } \right)^{\al/Q} < C \si || u||_{\infty}^{\al/Q} meas(B_{\la}) <   C \si || u||_{\infty}^{\al/Q} \var.
 \ee
 Hence
 \bea
 \int\limits_{A_\la} |\bu_{\la^2}(\zeta) - \bw_{\la^2}(\zeta)|^q \,d\zeta 
 &\leq &  C \si    ||u||_\infty^{\al/Q}  \var \int\limits_{A_\la} ( \bu_{\la^2}(\zeta) - {\bw_{\la^2}}(\zeta))^{q} \, d\zeta \nonumber \\
 &&  + C \si || u||_{\infty}^{\al/Q} \var ||( u - w) ||_{L^{q}(B_\la)} ||  (\bu_{\la^2}-  \bw_{\la^2})^{q-1}   ||_{L^l(A_\la)} \eea
 i.e.,
 \be \label{a6}
(1 -  C \si    ||u||_\infty^{\al/Q}  \var  ) \int\limits_{A_\la} |\bu_{\la^2}(\zeta) - \bw_{\la^2}(\zeta)|^q \,d\zeta  \leq C \si || u||_{\infty}^{\al/Q} \var ||( u - w) ||_{L^{q}(B_\la)} ||  (\bu_{\la^2}-  \bw_{\la^2})^{q-1}   ||_{L^l(A_\la)}. \ee
  Similarly,
  \bea \label{a7} 
  (1 -  C \si    ||u||_\infty^{\al/Q}  \var  ) \int\limits_{B_\la} |u(\zeta) - w(\zeta)|^q \,d\zeta  \leq C \si || u||_{\infty}^{\al/Q} \var ||( \bu_{\la^2} - \bw_{\la^2}) ||_{L^{q}(A_\la)} ||  (u-  w)^{q-1}   ||_{L^l(B_\la)}.\nonumber\\
\eea  
     Multiplying (\ref{a6}) and (\ref{a7}) we get
    \bea
 &&(1 -  C \si    ||u||_\infty^{\al/Q}  \var  )^2 \left( \int\limits_{A_\la} |\bu_{\la^2}(\zeta) - \bw_{\la^2}(\zeta)|^q \,d\zeta \right) \left(\int\limits_{B_\la} |u(\zeta) - w(\zeta)|^q \,d\zeta  \right)\nonumber\\
 &\leq &  (C \si || u||_{\infty}^{\al/Q} \var)^2 
 ||( u - w) ||_{L^{q}(B_\la)} ||  (\bu_{\la^2}-  \bw_{\la^2})^{q-1}   ||_{L^l(A_\la)}
 ||( \bu_{\la^2} - \bw_{\la^2}) ||_{L^{q}(A_\la)} ||  (u-  w)^{q-1}   ||_{L^l(B_\la)} \nonumber \\
 & = &  (C \si || u||_{\infty}^{\al/Q} \var)^2 \left( \int\limits_{A_\la} |\bu_{\la^2}(\zeta) - \bw_{\la^2}(\zeta)|^q \,d\zeta \right) \left(\int\limits_{B_\la} |u(\zeta) - w(\zeta)|^q \,d\zeta  \right) \eea
and we  get a contradiction to the definition of $\La$ by choosing $\var$ sufficiently small. 
Hence $\La = 1$ which gives that 
\be meas \{ \xi \in B(0,\la) :  u(\xi ) \geq w(\xi)  \} = 0 \mbox{ for all } 0 \leq \la \leq 1. \ee
Arguing as in (\ref{-4})-(\ref{-5}) and using (\ref{contra}), it follows that  \be \label{a4}
u(\xi) < w(\xi) \mbox{ and } \bu = \bu_{\la^2 =1} < \bw_{\la^2=1} = \bw \quad \mbox{  for all }\xi \in B(0,1). \ee But $u(\xi) = \bu(\xi)$ and $ w(\xi) =\frac{s_0^{Q-\al}}{\rho^{Q-\al}} w (s_0^2\hat \eta)$  and the fact that both $u$ and $w$ are solutions of the integral equation (\ref{m}) implies that 
\bea  
u_\infty = \bu(0) &=& \int\limits_{\HH^n} G_\al(0, \xi) \bu(\xi)^{\frac{Q+\al}{Q-\al}} \, d\xi  \nonumber \\
&=& \int\limits_{\HH^n} G_\al(0, \xi) u(\xi)^{\frac{Q+\al}{Q-\al}} \, d\xi \nonumber \\
& < & \int\limits_{\HH^n} G_\al(0, \xi)  w(\xi)^{\frac{Q+\al}{Q-\al}} \, d\xi \nonumber \\
& = & \int\limits_{\HH^n} G_\al(0, \xi) [\frac{s_0^{Q-\al}}{\rho^{Q-\al}} w (s_0^2\hat \xi)]^{\frac{Q+\al}{Q-\al}} \, d\xi = \int\limits_{\HH^n} G_\al(\zeta, \xi)  w_{s_0}( \xi)^{\frac{Q+\al}{Q-\al}} \, d\xi \nonumber \\
& = & \bw_{s_0^2}(0) = w_\infty \eea
a contradiction! Hence \be w(0) \leq u(0). \ee Interchanging the roles of $u$ and $w$ and repeating the proof, we conclude that $u(0) \leq w(0)$ and hence 
\be\label{a5} u(0) = w(0). \ee
In particular, (\ref{a5}) implies that $w$ is also inversion symmetric with respect to the unit CC-sphere i.e., 
\be w (\xi) = \bw(\xi) \quad \mbox{ for all } \xi \in \HH^n. \ee

The above proof can be repeated to  further conclude that   
\be  u \equiv w \mbox{ in } \HH^n. \ee
\qed

{\bf Proof of Theorem \ref{uniq}:\/} Let \be 
u_0 = C_0 ( t^2 + (r^2 + 1)^2)^{- \frac{Q-\al}{4} } = C_0 |\om +i|^{\frac{Q-\al}{2}}\ee
denote the standard solution of (\ref{m}) centered at the origin. It can be directly verified that $u_0 (\xi) = \bu_0
(\xi) $ for all $\xi \in \HH^n$. 
If $w$ is  any other solution of the equation (\ref{m}), then from Theorem \ref{symm} it follows 
\be w_\infty := \lim\limits_{\rho \to \infty} \rho^{Q-\al} w(\xi) \mbox{ exists } \ee
and $w$ is $\HH$-symmetric with respect to some hyperplane ${\mathcal H}_{\xi_0}$. \\
Now consider $v(\xi) = w_s:=s^{\frac{Q-\al}{2}} w (sr, s^2 t)$. Then 
$v_\infty = \lim\limits_{\rho \to \infty} \rho^{Q-\al} v(\xi) = \frac{w_\infty}{s^{\frac{Q-\al}{2}}}$. \\
Choose $s^{\frac{Q-\al}{2}} = \frac{w_\infty}{(u_0)_\infty}$ so that we get
\[v_\infty = (u_0)_\infty= u_0(0).\] 
Again, from Theorem \ref{symm} we know that after a translation by  $\xi_0$, $v_1(\xi) = v(\xi_0\xi )$  is cylindrical i.e., $v_1(\xi) = v_1(r, t)$ and $(v_1)_\infty = (u_0)_\infty$.  Hence, from Lemma \ref{uniq}, 
$ v_1 \equiv u_0 $, i.e., 
\[ v(\xi) = u_0( \xi_0^{-1} \xi) \mbox{ or }  s^{\frac{Q-\al}{2}} w (sr, s^2 t) = u_0( \xi_0^{-1} \xi). \]
It follows that
\[w(\xi) =  s^{-\frac{Q-\al}{2}} u_0( \xi_0^{-1}(\delta_{1/s}\xi)) \]
which proves the uniqueness. 

\qed 


{\bf Acknowledgements:\/} The first author would like to dedicate this paper to her teachers and mentors, particularly late Professor Abbas Bahri who had always encouraged and supported  her efforts to understand the Heisenberg group. She also acknowledges the International Center for Theoretical Physics ( ICTP, Trieste) for granting Associateship during the period 2018-2023, facilitating her research  during her visits to the ICTP. The first author thanks Claudio Afeltra for pointing out the error in the initial version of the paper.

 \end{document}